# COUPLED PARAXIAL WAVE EQUATIONS IN RANDOM MEDIA IN THE WHITE-NOISE REGIME[1]

By Josselin Garnier and Knut Sølna[2]

*Université Paris VII and University of California at Irvine*

In this paper the reflection and transmission of waves by a three-dimensional random medium are studied in a white-noise and paraxial regime. The limit system derives from the acoustic wave equations and is described by a coupled system of random Schrödinger equations driven by a Brownian field whose covariance is determined by the two-point statistics of the fluctuations of the random medium. For the reflected and transmitted fields the associated Wigner distributions and the autocorrelation functions are determined by a closed system of transport equations. The Wigner distribution is then used to describe the enhanced backscattering phenomenon for the reflected field.

**1. Introduction.** The paraxial wave equation, in homogeneous or in random media, is a model used for many applications, for instance in communication and imaging [19]. It has the form of an evolution equation that describes waves propagating along a privileged axis and it can be obtained by neglecting backscattering. Its simplicity, compared to the full three-dimensional wave equation, enables analysis of many important phenomena, such as laser beam propagation [9, 12, 16, 23], time reversal in random media [5, 11], underwater acoustics [24] or migration problems in geophysics [6].

The derivation of the paraxial model in homogeneous media is well understood, and it can also be justified in heterogeneous media for small variations of the wave speed [2, 4]. However, it is not clear whether the paraxial (parabolic) approximation is still valid in a scaling regime in which the medium fluctuations are rapid and can be approximated by a white-noise

Received December 2007; revised April 2008.
[1]Supported by ONR Grant N00014-02-1-0089 and DARPA Grant N00014-05-1-0442.
[2]Supported by NSF Grant DMS-03-07011 and the Sloan Foundation.

term. The main motivations for studying the white-noise paraxial wave equation are (i) it appears as a very natural model in many applications where the correlation length of the medium is relatively small, in particular much smaller than the propagation distance, (ii) it allows for the use of Itô's stochastic calculus, which in turn enables the closure of the hierarchy of moment equations [15] and thereby analysis of important wave propagation problems, such as the star scintillation due to atmospheric turbulence [25].

If the paraxial approximation and the white-noise approximation can be justified simultaneously, then the conjecture [21, 23] is that the limit equation should take the form of the random Schrödinger equation studied in particular in [8]. The proof of the convergence of the solution of the wave equation in random media to the solution of the white-noise paraxial equation was obtained in the case of stratified weakly fluctuating media in [1]. In our paper we consider the transmission and reflection of acoustic waves by a slab of medium whose parameters have three-dimensional random fluctuations and whose end is either transparent or a strong interface. This model is particularly interesting in the context of optical coherence tomography [26]. We analyze a wave propagation regime in which the paraxial and white-noise approximations are valid. In this regime we obtain a system of coupled Schrödinger equations driven by a Brownian field that fully determines the statistics of the transmitted and reflected waves. As a corollary we compute explicitly the two-point statistics of the transmitted and reflected waves. These results show that the often used "independent approach" for the reflected wave (in which the statistics of the forward- and backward-propagating waves are assumed to be independent [26]) is valid if and only if the transverse correlation radius of the fluctuations of the random medium is smaller than the initial beam width. Finally, we use the coupled Schrödinger equations to give a rigorous account for the enhanced backscattering or weak localization phenomenon [27] and we compute explicitly the enhancement factor and the shape of the enhanced backscattering cone.

**2. The transmission and reflection operators.** We consider linear acoustic waves propagating in $1 + d$ spatial dimensions with random medium fluctuations. The governing equations are

$$(2.1) \qquad \rho(z, \mathbf{x}) \frac{\partial \mathbf{u}}{\partial t} + \nabla p = \mathbf{F}, \qquad \frac{1}{K(z, \mathbf{x})} \frac{\partial p}{\partial t} + \nabla \cdot \mathbf{u} = 0,$$

where $p$ is the pressure field, $\mathbf{u}$ is the velocity field, $\rho$ is the density of the medium, $K$ is the bulk modulus of the medium and $(z, \mathbf{x}) \in \mathbb{R} \times \mathbb{R}^d$ are the space coordinates. The source is modelled by the forcing term $\mathbf{F}$. We consider in this paper the situation in which a random slab occupying the section $z \in (0, L)$ is sandwiched in between two homogeneous half-spaces. The source, $\mathbf{F}$, is located outside of the slab at $z = z_0$, $z_0 > L$. We shall refer to waves propagating in the direction with a positive $z$ component as right-propagating waves. The medium fluctuations in the random slab $(0, L)$ vary rapidly in space while the "background" medium is constant. The medium is assumed to be matched at the right boundary $z = L$. We consider a possible mismatch at the boundary $z = 0$ and denote the background medium parameters by $\rho_0$ and $K_0$ in the half-space $z \leq 0$ and by $\rho_1$ and $K_1$ in the half-space $z > 0$:

$$\frac{1}{K(z,\mathbf{x})} = \begin{cases} K_0^{-1}, & \text{if } z \leq 0, \\ K_1^{-1}(1 + \nu_K(z, \mathbf{x})), & \text{if } z \in (0, L), \\ K_1^{-1}, & \text{if } z \geq L, \end{cases}$$

$$\rho(z, \mathbf{x}) = \begin{cases} \rho_0, & \text{if } z \leq 0, \\ \rho_1, & \text{if } z \in (0, L), \\ \rho_1, & \text{if } z \geq L, \end{cases}$$

where the random field $\nu_K(z, \mathbf{x})$ models the medium fluctuations, whose correlation length is $l_K$. The source has the form

$$\mathbf{F}(t, z, \mathbf{x}) = f_s(t, \mathbf{x})\delta(z - z_0)\mathbf{e}_z,$$

where $\mathbf{e}_z$ is the unit vector pointing in the $z$-direction, $z_0 > L$ is the source position. We denote by $\omega_0$ the typical frequency of the source term $f_s$ and by $R_0$ the diameter of its spatial support (which gives the initial beam width). The typical wavelength associated with the typical frequency $\omega_0$ is $\lambda_0 = 2\pi c_1/\omega_0$, for $c_1 = \sqrt{K_1/\rho_1}$ the background speed for $z > 0$, which is of the same order as the background speed $c_0 = \sqrt{K_0/\rho_0}$ in the half-space $z \leq 0$.

We can now introduce the scaling regime that we consider in this paper:

(1) We assume that the correlation length $l_K$ of the medium is much smaller than the propagation distance $L$. We denote by $\varepsilon^2$ the ratio between the correlation length and the typical propagation distance.

(2) We assume that the transverse width of the source $R_0$ and the correlation length of the medium $l_K$ are of the same order. This means that we assume that the ratio $R_0/L$ is of order $\varepsilon^2$. This scaling is motivated by the fact that, in this regime, there is a nontrivial interaction between the fluctuations of the medium and the beam.

(3) We assume that the typical wavelength $\lambda_0$ is much smaller than the propagation distance $L$; more precisely, we assume that the ratio $\lambda_0/L$ is of order $\varepsilon^4$. This high-frequency scaling is motivated by the following considerations. The Rayleigh length for a beam with initial width $R_0$ and central wavelength $\lambda_0$ is of the order of $R_0^2/\lambda_0$ in absence of random fluctuations (the Rayleigh length is the distance from beam waist where the beam area is doubled by diffraction). In order to get a Rayleigh length of the order of the propagation distance $L$, the ratio $\lambda_0/L$ must be of order $\varepsilon^4$ since $R_0/L \sim \varepsilon^2$.

Henceforth we shall assume nondimensionalized units chosen such that the background bulk modulus $K_1$ and density $\rho_1$ in the half-space $z > 0$ are 1, hence, the background speed $c_1 = \sqrt{K_1/\rho_1}$ and impedance $Z_1 = \sqrt{K_1 \rho_1}$ are also equal to 1. If we consider the propagation distance, $L$, as our reference distance of order 1 in this scaled regime, then the source has the form

$$(2.2) \qquad \mathbf{F}(t, z, \mathbf{x}) = f\left(\frac{t}{\varepsilon^4}, \frac{\mathbf{x}}{\varepsilon^2}\right) \delta(z - z_0) \mathbf{e}_z,$$

where $f(t, \mathbf{x})$ is the normalized source shape function (with time and spatial scales of variations of order 1), and the medium fluctuations have the form

$$\frac{1}{K(z, \mathbf{x})} = \begin{cases} K_0^{-1}, & \text{if } z \leq 0, \\ 1 + \varepsilon^3 \nu\left(\frac{z}{\varepsilon^2}, \frac{\mathbf{x}}{\varepsilon^2}\right), & \text{if } z \in (0, L), \\ 1, & \text{if } z \geq L, \end{cases}$$

$$\rho(z, \mathbf{x}) = \begin{cases} \rho_0, & \text{if } z \leq 0, \\ 1, & \text{if } z \in (0, L), \\ 1, & \text{if } z \geq L, \end{cases}$$

where the zero-mean, stationary random field $\nu$ has a correlation length of order 1 and standard deviation of order 1. We also assume that it satisfies strong mixing conditions in $z$. Here the amplitude $\varepsilon^3$ of the fluctuations has been chosen so as to obtain an effective regime of order 1 when $\varepsilon$ goes to zero. That is, if the magnitude of the fluctuations is smaller than $\varepsilon^3$, then the wave would propagate as if the medium were homogeneous, while if the order of magnitude is larger, then the wave would not penetrate the slab. The scaling that we consider here corresponds to the physically most interesting situation.

Since both the medium and the source have transverse spatial variations at the scale $\varepsilon^2$, it is convenient to rescale the transverse variable $\mathbf{x}/\varepsilon^2 \to \mathbf{x}$ and to introduce the rescaled fields $\mathbf{u}^\varepsilon$ and $p^\varepsilon$:

$$(2.3) \qquad \mathbf{u}^\varepsilon(t, z, \mathbf{x}) = \mathbf{u}(t, z, \varepsilon^2 \mathbf{x}), \qquad p^\varepsilon(t, z, \mathbf{x}) = p(t, z, \varepsilon^2 \mathbf{x}).$$

The reader should keep in mind that thus, in the discussion below, when we refer to the transversal spatial parameter $\mathbf{x}$ it corresponds to $\varepsilon^2 \mathbf{x}$ in the original coordinates. The rescaled fields satisfy in the region $z \in (-\infty, 0]$:

$$(2.4) \qquad \rho_0 \frac{\partial \mathbf{u}^\varepsilon}{\partial t} + \begin{bmatrix} \partial_z \\ \varepsilon^{-2} \nabla_\mathbf{x} \end{bmatrix} p^\varepsilon = \mathbf{0}, \qquad \frac{1}{K_0} \frac{\partial p^\varepsilon}{\partial t} + \begin{bmatrix} \partial_z \\ \varepsilon^{-2} \nabla_\mathbf{x} \end{bmatrix} \cdot \mathbf{u}^\varepsilon = 0,$$

where $\nabla_\mathbf{x}$ stands for the gradient with respect to the transverse spatial variables $\mathbf{x}$. In the random slab $z \in (0, L)$ the fields satisfy

$$(2.5) \qquad \begin{aligned} \frac{\partial \mathbf{u}^\varepsilon}{\partial t} + \begin{bmatrix} \partial_z \\ \varepsilon^{-2} \nabla_\mathbf{x} \end{bmatrix} p^\varepsilon &= \mathbf{0}, \\ \left(1 + \varepsilon^3 \nu\left(\frac{z}{\varepsilon^2}, \mathbf{x}\right)\right) \frac{\partial p^\varepsilon}{\partial t} + \begin{bmatrix} \partial_z \\ \varepsilon^{-2} \nabla_\mathbf{x} \end{bmatrix} \cdot \mathbf{u}^\varepsilon &= 0, \end{aligned}$$

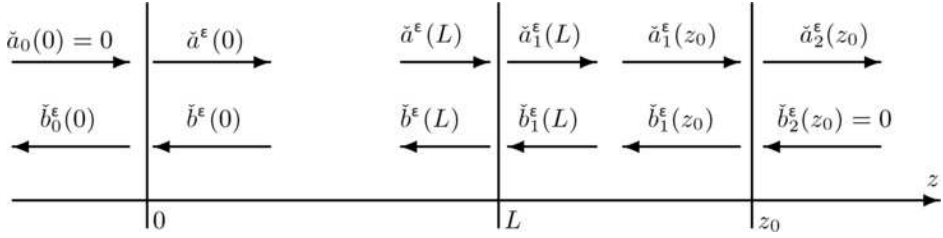

FIG. 1. *Boundary conditions for the modes in the presence of an interface at $z = 0$, a random slab $(0, L)$ and a source at $z = z_0$.*

and in the region $z \in [L, \infty)$ (in which the source is located):

$$\frac{\partial \mathbf{u}^\varepsilon}{\partial t} + \begin{bmatrix} \partial_z \\ \varepsilon^{-2} \nabla_\mathbf{x} \end{bmatrix} p^\varepsilon = f\left(\frac{t}{\varepsilon^4}, \mathbf{x}\right) \delta(z - z_0) \begin{bmatrix} 1 \\ \mathbf{0} \end{bmatrix}, \qquad \frac{\partial p^\varepsilon}{\partial t} + \begin{bmatrix} \partial_z \\ \varepsilon^{-2} \nabla_\mathbf{x} \end{bmatrix} \cdot \mathbf{u}^\varepsilon = 0.$$

In the Fourier domain (with respect to time), the first equations of these systems give the expressions of the velocity fields in terms of the pressure fields.

We first consider the wave field in the homogeneous half-space $z \leq 0$, which allows us to introduce the standard paraxial wave equation in homogeneous medium. The wave field satisfies (2.4) with the wave speed $c_0 = \sqrt{K_0/\rho_0}$. Following [18] (in a different scaling) we introduce the complex amplitudes $\check{a}_0^\varepsilon$ and $\check{b}_0^\varepsilon$ of the right- and left-propagating modes (see Figure 1):

$$\check{a}_0^\varepsilon(k, z, \mathbf{x}) = \frac{c_0}{2} \left[ \int \left( \frac{1}{\varepsilon^4} p^\varepsilon(t, z, \mathbf{x}) + \frac{1}{ik} \frac{\partial p^\varepsilon}{\partial z}(t, z, \mathbf{x}) \right) e^{ic_0 kt/\varepsilon^4} \, dt \right] e^{-ikz/\varepsilon^4},$$

$$\check{b}_0^\varepsilon(k, z, \mathbf{x}) = \frac{c_0}{2} \left[ \int \left( \frac{1}{\varepsilon^4} p^\varepsilon(t, z, \mathbf{x}) - \frac{1}{ik} \frac{\partial p^\varepsilon}{\partial z}(t, z, \mathbf{x}) \right) e^{ic_0 kt/\varepsilon^4} \, dt \right] e^{ikz/\varepsilon^4}.$$

They are such that the pressure field in the region $z \leq 0$ can be written as

$$p^\varepsilon(t, z, \mathbf{x}) = \frac{1}{2\pi} \int (\check{a}_0^\varepsilon(k, z, \mathbf{x}) e^{ikz/\varepsilon^4} + \check{b}_0^\varepsilon(k, z, \mathbf{x}) e^{-ikz/\varepsilon^4}) e^{-ic_0 kt/\varepsilon^4} \, dk,$$

and they satisfy

$$\frac{\partial \check{a}_0^\varepsilon}{\partial z}(k, z, \mathbf{x}) e^{ikz/\varepsilon^4} + \frac{\partial \check{b}_0^\varepsilon}{\partial z}(k, z, \mathbf{x}) e^{-ikz/\varepsilon^4} = 0.$$

Using (2.4), we find that they also satisfy the coupled mode equations

$$\frac{\partial \check{a}_0^\varepsilon}{\partial z} = \frac{i}{2k} \Delta_\mathbf{x} \check{a}_0^\varepsilon + e^{-2ikz/\varepsilon^4} \frac{i}{2k} \Delta_\mathbf{x} \check{b}_0^\varepsilon, \qquad \frac{\partial \check{b}_0^\varepsilon}{\partial z} = -e^{2ikz/\varepsilon^4} \frac{i}{2k} \Delta_\mathbf{x} \check{a}_0^\varepsilon - \frac{i}{2k} \Delta_\mathbf{x} \check{b}_0^\varepsilon,$$

where $\Delta_\mathbf{x}$ is the transverse Laplacian. In the limit $\varepsilon \to 0$, the cross terms (proportional to $e^{\pm 2ikz/\varepsilon^4}$) average out to zero and we get the two uncoupled

paraxial wave equations

$$\frac{\partial \check{a}_0^\varepsilon}{\partial z} = \frac{i}{2k}\Delta_{\mathbf{x}}\check{a}_0^\varepsilon, \qquad \frac{\partial \check{b}_0^\varepsilon}{\partial z} = -\frac{i}{2k}\Delta_{\mathbf{x}}\check{b}_0^\varepsilon.$$

Taking into account the fact that there is no source in the half-space $z \leq 0$, and therefore no right-going wave, we obtain the general expression of the wave in the left homogeneous half-space

$$(2.6) \qquad p^\varepsilon(t,z,\mathbf{x}) = \frac{1}{2\pi}\int \check{b}_0^\varepsilon(k,z,\mathbf{x})e^{-ikz/\varepsilon^4}e^{-ic_0kt/\varepsilon^4}\,dk, \qquad z \leq 0.$$

Similarly, the wave fields in the homogeneous regions $(L,z_0)$ and $(z_0,\infty)$ have the forms

$$p^\varepsilon(t,z,\mathbf{x})$$
$$= \begin{cases} \dfrac{1}{2\pi}\int (\check{a}_1^\varepsilon(k,z,\mathbf{x})e^{ikz/\varepsilon^4} + \check{b}_1^\varepsilon(k,z,\mathbf{x})e^{-ikz/\varepsilon^4})e^{-ikt/\varepsilon^4}\,dk, & z \in (L,z_0), \\ \dfrac{1}{2\pi}\int \check{a}_2^\varepsilon(k,z,\mathbf{x})e^{ikz/\varepsilon^4}e^{-ikt/\varepsilon^4}\,dk, & z > z_0, \end{cases}$$

respectively. Here we used the fact that there is no source and therefore no left-going wave in the region $z > z_0$; see Figure 1. We can also use the jump conditions across the source position $z = z_0$ to obtain the relations

$$(2.7) \qquad \check{b}_1^\varepsilon(k,z_0,\mathbf{x}) = -\tfrac{1}{2}e^{ikz_0/\varepsilon^4}\check{f}(k,\mathbf{x}),$$

$$(2.8) \qquad \check{a}_2^\varepsilon(k,z_0,\mathbf{x}) - \check{a}_1^\varepsilon(k,z_0,\mathbf{x}) = \tfrac{1}{2}e^{-ikz_0/\varepsilon^4}\check{f}(k,\mathbf{x}).$$

By solving the paraxial wave equation for $\check{b}_1^\varepsilon$, we obtain the expression for the complex amplitude of the wave impinging on the random slab at $z = L$:

$$(2.9) \qquad \check{b}_1^\varepsilon(k,L,\mathbf{x}) = e^{ikz_0/\varepsilon^4}\check{b}_{\mathrm{inc}}(k,\mathbf{x}),$$

$$(2.10) \qquad \check{b}_{\mathrm{inc}}(k,\mathbf{x}) = -\frac{1}{2(2\pi)^d}\int \hat{f}(k,\boldsymbol{\kappa})e^{i|\boldsymbol{\kappa}|^2(L-z_0)/(2k)+i\boldsymbol{\kappa}\cdot\mathbf{x}}\,d\boldsymbol{\kappa},$$

where the Fourier transforms are defined by

$$(2.11) \quad \check{f}(k,\mathbf{x}) = \int f(t,\mathbf{x})e^{ikt}\,dt, \qquad \hat{f}(k,\boldsymbol{\kappa}) = \int \check{f}(k,\mathbf{x})e^{-i\boldsymbol{\kappa}\cdot\mathbf{x}}\,d\mathbf{x}.$$

The pressure field in the region $z \in (0,L)$ can be written as

$$p^\varepsilon(t,z,\mathbf{x}) = \frac{1}{2\pi}\int (\check{a}^\varepsilon(k,z,\mathbf{x})e^{ikz/\varepsilon^4} + \check{b}^\varepsilon(k,z,\mathbf{x})e^{-ikz/\varepsilon^4})e^{-ikt/\varepsilon^4}\,dk,$$

with the complex amplitudes $\check{a}^\varepsilon$ and $\check{b}^\varepsilon$ of the right- and left-propagating modes given explicitly by

$$\check{a}^\varepsilon(k,z,\mathbf{x}) = \frac{1}{2}\left[\int \left(\frac{1}{\varepsilon^4}p^\varepsilon(t,z,\mathbf{x}) + \frac{1}{ik}\frac{\partial p^\varepsilon}{\partial z}(t,z,\mathbf{x})\right)e^{ikt/\varepsilon^4}\,dt\right]e^{-ikz/\varepsilon^4},$$

$$\check{b}^\varepsilon(k,z,\mathbf{x}) = \frac{1}{2}\left[\int \left(\frac{1}{\varepsilon^4}p^\varepsilon(t,z,\mathbf{x}) - \frac{1}{ik}\frac{\partial p^\varepsilon}{\partial z}(t,z,\mathbf{x})\right)e^{ikt/\varepsilon^4}\,dt\right]e^{ikz/\varepsilon^4}.$$

Using (2.5) we obtain the following mode coupling equations:

$$\frac{\partial \check{a}^\varepsilon}{\partial z} = \left(\frac{ik}{2\varepsilon}\nu\left(\frac{z}{\varepsilon^2}, \mathbf{x}\right) + \frac{i}{2k}\Delta_\mathbf{x}\right)\check{a}^\varepsilon$$
$$+ e^{-2ikz/\varepsilon^4}\left(\frac{ik}{2\varepsilon}\nu\left(\frac{z}{\varepsilon^2}, \mathbf{x}\right) + \frac{i}{2k}\Delta_\mathbf{x}\right)\check{b}^\varepsilon, \tag{2.12}$$

$$\frac{\partial \check{b}^\varepsilon}{\partial z} = -e^{2ikz/\varepsilon^4}\left(\frac{ik}{2\varepsilon}\nu\left(\frac{z}{\varepsilon^2}, \mathbf{x}\right) + \frac{i}{2k}\Delta_\mathbf{x}\right)\check{a}^\varepsilon$$
$$- \left(\frac{ik}{2\varepsilon}\nu\left(\frac{z}{\varepsilon^2}, \mathbf{x}\right) + \frac{i}{2k}\Delta_\mathbf{x}\right)\check{b}^\varepsilon. \tag{2.13}$$

This system is valid in $z \in (0, L)$ and it is complemented with the following boundary conditions at $z = 0$ and $z = L$:

$$\check{b}^\varepsilon(k, z = L, \mathbf{x}) = e^{ikz_0/\varepsilon^4}\check{b}_{\mathrm{inc}}(k, \mathbf{x}), \tag{2.14}$$

$$\check{a}^\varepsilon(k, z = 0, \mathbf{x}) = \mathcal{R}_0 \check{b}^\varepsilon(k, z = 0, \mathbf{x}), \tag{2.15}$$

where $\mathcal{R}_0 = (Z_0 - 1)/(Z_0 + 1)$ is the reflection coefficient of the interface at $z = 0$ and $Z_0 = \sqrt{K_0\rho_0}$ is the impedance of the left homogeneous half-space. These boundary conditions are obtained from the continuity relations of the fields $p^\varepsilon$ and $\mathbf{e}_z \cdot \mathbf{u}^\varepsilon$ at $z = 0$ and $z = L$. The continuity relations also give the expressions for the complex amplitudes of the transmitted field $\check{b}_0^\varepsilon$ in the region $z \leq 0$ and of the reflected field $\check{a}_1^\varepsilon$ in the region $z \geq L$:

$$\check{b}_0^\varepsilon(k, z = 0, \mathbf{x}) = \mathcal{T}_0 \check{b}^\varepsilon(k, z = 0, \mathbf{x}), \tag{2.16}$$

$$\check{a}_1^\varepsilon(k, z = L, \mathbf{x}) = \check{a}^\varepsilon(k, z = L, \mathbf{x}), \tag{2.17}$$

where $\mathcal{T}_0 = 2Z_0^{1/2}/(1 + Z_0)$ is the transmission coefficient of the interface at $z = 0$. If there is no impedance contrast $Z_0 = 1$, then $\mathcal{T}_0 = 1$ and $\mathcal{R}_0 = 0$ and the boundary condition (2.15) reads $\check{a}^\varepsilon(k, z = 0, \mathbf{x}) = 0$. This is the radiation condition expressing the fact that there is no wave incoming from $-\infty$.

We now make use of an invariant imbedding step and introduce transmission and reflection operators. First, we define the lateral Fourier modes

$$\hat{a}^\varepsilon(k, z, \boldsymbol{\kappa}) = \int \check{a}^\varepsilon(k, z, \mathbf{x}) e^{-i\boldsymbol{\kappa}\cdot\mathbf{x}}\, d\mathbf{x},$$
$$\hat{b}^\varepsilon(k, z, \boldsymbol{\kappa}) = \int \check{b}^\varepsilon(k, z, \mathbf{x}) e^{-i\boldsymbol{\kappa}\cdot\mathbf{x}}\, d\mathbf{x}, \tag{2.18}$$

and make the ansatz

$$\hat{b}_0^\varepsilon(k, 0, \boldsymbol{\kappa}) = \int \widehat{\mathcal{T}}^\varepsilon(k, z, \boldsymbol{\kappa}, \boldsymbol{\kappa}')\hat{b}^\varepsilon(k, z, \boldsymbol{\kappa}')\, d\boldsymbol{\kappa}',$$
$$\hat{a}^\varepsilon(k, z, \boldsymbol{\kappa}) = \int \widehat{\mathcal{R}}^\varepsilon(k, z, \boldsymbol{\kappa}, \boldsymbol{\kappa}')\hat{b}^\varepsilon(k, z, \boldsymbol{\kappa}')\, d\boldsymbol{\kappa}'. \tag{2.19}$$

Using the mode coupling equations (2.12)–(2.13) we find that the operators $\widehat{\mathcal{T}}^\varepsilon$ and $\widehat{\mathcal{R}}^\varepsilon$ satisfy

$$
\begin{aligned}
(2.20) \quad \frac{d}{dz}&\widehat{\mathcal{R}}^\varepsilon(k,z,\boldsymbol{\kappa},\boldsymbol{\kappa}') \\
&= e^{-2ikz/\varepsilon^4}\widehat{\mathcal{L}}^\varepsilon(k,z,\boldsymbol{\kappa},\boldsymbol{\kappa}') \\
&\quad + e^{2ikz/\varepsilon^4}\iint \widehat{\mathcal{R}}^\varepsilon(k,z,\boldsymbol{\kappa},\boldsymbol{\kappa}_1)\widehat{\mathcal{L}}^\varepsilon(k,z,\boldsymbol{\kappa}_1,\boldsymbol{\kappa}_2) \\
&\qquad \times \widehat{\mathcal{R}}^\varepsilon(k,z,\boldsymbol{\kappa}_2,\boldsymbol{\kappa}')\,d\boldsymbol{\kappa}_1\,d\boldsymbol{\kappa}_2 \\
&\quad + \int \widehat{\mathcal{L}}^\varepsilon(k,z,\boldsymbol{\kappa},\boldsymbol{\kappa}_1)\widehat{\mathcal{R}}^\varepsilon(k,z,\boldsymbol{\kappa}_1,\boldsymbol{\kappa}')\,d\boldsymbol{\kappa}_1 \\
&\quad + \int \widehat{\mathcal{R}}^\varepsilon(k,z,\boldsymbol{\kappa},\boldsymbol{\kappa}_1)\widehat{\mathcal{L}}^\varepsilon(k,z,\boldsymbol{\kappa}_1,\boldsymbol{\kappa}')\,d\boldsymbol{\kappa}_1,
\end{aligned}
$$

$$
\begin{aligned}
(2.21) \quad \frac{d}{dz}&\widehat{\mathcal{T}}^\varepsilon(k,z,\boldsymbol{\kappa},\boldsymbol{\kappa}') \\
&= \int \widehat{\mathcal{T}}^\varepsilon(k,z,\boldsymbol{\kappa},\boldsymbol{\kappa}_1)\widehat{\mathcal{L}}^\varepsilon(k,z,\boldsymbol{\kappa}_1,\boldsymbol{\kappa}')\,d\boldsymbol{\kappa}_1 \\
&\quad + e^{2ikz/\varepsilon^4}\iint \widehat{\mathcal{T}}^\varepsilon(k,z,\boldsymbol{\kappa},\boldsymbol{\kappa}_1)\widehat{\mathcal{L}}^\varepsilon(k,z,\boldsymbol{\kappa}_1,\boldsymbol{\kappa}_2) \\
&\qquad \times \widehat{\mathcal{R}}^\varepsilon(k,z,\boldsymbol{\kappa}_2,\boldsymbol{\kappa}')\,d\boldsymbol{\kappa}_1\,d\boldsymbol{\kappa}_2,
\end{aligned}
$$

where we have defined

$$
(2.22) \quad \widehat{\mathcal{L}}^\varepsilon(k,z,\boldsymbol{\kappa}_1,\boldsymbol{\kappa}_2) = -\frac{i}{2k}|\boldsymbol{\kappa}_1|^2 \delta(\boldsymbol{\kappa}_1-\boldsymbol{\kappa}_2) + \frac{ik}{2(2\pi)^d\varepsilon}\hat{\nu}\Big(\frac{z}{\varepsilon^2},\boldsymbol{\kappa}_1-\boldsymbol{\kappa}_2\Big),
$$

with $\hat{\nu}(z,\boldsymbol{\kappa})$ the partial Fourier transform (in **x**) of $\nu(z,\mathbf{x})$. This system is complemented with the initial conditions at $z=0$, which are obtained from (2.15) and (2.16):

$$
\widehat{\mathcal{R}}^\varepsilon(k,z=0,\boldsymbol{\kappa},\boldsymbol{\kappa}') = \mathcal{R}_0\delta(\boldsymbol{\kappa}-\boldsymbol{\kappa}'), \qquad \widehat{\mathcal{T}}^\varepsilon(k,z=0,\boldsymbol{\kappa},\boldsymbol{\kappa}') = \mathcal{T}_0\delta(\boldsymbol{\kappa}-\boldsymbol{\kappa}').
$$

The transmission and reflection operators evaluated at $z=L$ carry all the relevant information about the random medium from the point of view of the transmitted and reflected waves, which are our main quantities of interest.

Our objective in the next sections is to characterize the transmitted wave field

$$
\begin{aligned}
(2.23) \quad p_{\mathrm{tr}}^\varepsilon(s,\mathbf{x}) &= p^\varepsilon(z_0+\varepsilon^4 s, z=0,\mathbf{x}) \\
&= \frac{1}{(2\pi)^{d+1}}\iiint \widehat{\mathcal{T}}^\varepsilon(k,L,\boldsymbol{\kappa},\boldsymbol{\kappa}')\hat{b}_{\mathrm{inc}}(k,\boldsymbol{\kappa}')\,d\boldsymbol{\kappa}'\,e^{i(\boldsymbol{\kappa}\cdot\mathbf{x}-ks)}\,d\boldsymbol{\kappa}\,dk,
\end{aligned}
$$

and the reflected wave field

$$p_{\text{ref}}^\varepsilon(s,\mathbf{x}) = p^\varepsilon(z_0 + L + \varepsilon^4 s, z = L, \mathbf{x})$$
$$(2.24) \qquad = \frac{1}{(2\pi)^{d+1}} \iiint \widehat{\mathcal{R}}^\varepsilon(k, L, \boldsymbol{\kappa}, \boldsymbol{\kappa}')\hat{b}_{\text{inc}}(k,\boldsymbol{\kappa}')\,d\boldsymbol{\kappa}'\,e^{i(\boldsymbol{\kappa}\cdot\mathbf{x}-ks)}\,d\boldsymbol{\kappa}\,dk.$$

Note that the wave field "fronts" are observed on the time scale of the source and around their respective expected arrival times ($z_0$ for the transmitted wave, and $L + z_0$ for the reflected wave, which corresponds to the travel time to go from the source at $z = z_0$ to the interface at $z = 0$ and back at the surface at $z = L$).

## 3. The random Schrödinger model.

3.1. *Statement of the main result.* We consider the transmitted and reflected fields $p_{\text{tr}}^\varepsilon$ and $p_{\text{ref}}^\varepsilon$ defined by (2.24)–(2.25) and use diffusion approximation theorems to identify a limit random Schrödinger model. The main result is the following one.

PROPOSITION 3.1. *The processes* $(p_{\text{tr}}^\varepsilon(s,\mathbf{x}), p_{\text{ref}}^\varepsilon(s,\mathbf{x}))_{s\in\mathbb{R},\mathbf{x}\in\mathbb{R}^d}$ *converge in distribution as* $\varepsilon \to 0$ *in the space* $C^0(\mathbb{R}, L^2(\mathbb{R}^d, \mathbb{R}^2)) \cap L^2(\mathbb{R}, L^2(\mathbb{R}^d, \mathbb{R}^2))$ *to the limit process* $(p_{\text{tr}}(s,\mathbf{x}), p_{\text{ref}}(s,\mathbf{x}))_{s\in\mathbb{R},\mathbf{x}\in\mathbb{R}^d}$

$$(3.1) \qquad p_{\text{tr}}(s,\mathbf{x}) = \frac{1}{2\pi}\iint \check{\mathcal{T}}(k, L, \mathbf{x}, \mathbf{x}')\check{b}_{\text{inc}}(k, \mathbf{x}')\,d\mathbf{x}'\,e^{-iks}\,dk,$$

$$(3.2) \qquad p_{\text{ref}}(s,\mathbf{x}) = \frac{1}{2\pi}\iint \check{\mathcal{R}}(k, L, \mathbf{x}, \mathbf{x}')\check{b}_{\text{inc}}(k, \mathbf{x}')\,d\mathbf{x}'\,e^{-iks}\,dk.$$

*Here* $C^0(\mathbb{R}, L^2(\mathbb{R}^d, \mathbb{R}^2))$ *is the space of continuous functions (in $s$) with values in* $L^2(\mathbb{R}^d, \mathbb{R}^2)$ *and* $L^2(\mathbb{R}, L^2(\mathbb{R}^d, \mathbb{R}^2)) = L^2(\mathbb{R} \times \mathbb{R}^d, \mathbb{R}^2)$. *The operators* $\check{\mathcal{T}}(k, z, \mathbf{x}, \mathbf{x}')$ *and* $\check{\mathcal{R}}(k, z, \mathbf{x}, \mathbf{x}')$ *are the solutions of the following Itô–Schrödinger diffusion models:*

$$(3.3)\ \ d\check{\mathcal{T}}(k, z, \mathbf{x}, \mathbf{x}') = \frac{i}{2k}\Delta_{\mathbf{x}'}\check{\mathcal{T}}(k, z, \mathbf{x}, \mathbf{x}')\,dz + \frac{ik}{2}\check{\mathcal{T}}(k, z, \mathbf{x}, \mathbf{x}') \circ dB(z, \mathbf{x}'),$$

$$(3.4) \quad \begin{aligned} d\check{\mathcal{R}}(k, z, \mathbf{x}, \mathbf{x}') &= \frac{i}{2k}(\Delta_{\mathbf{x}} + \Delta_{\mathbf{x}'})\check{\mathcal{R}}(k, z, \mathbf{x}, \mathbf{x}')\,dz \\ &\quad + \frac{ik}{2}\check{\mathcal{R}}(k, z, \mathbf{x}, \mathbf{x}') \circ (dB(z, \mathbf{x}) + dB(z, \mathbf{x}')), \end{aligned}$$

*with the initial conditions*

$$\check{\mathcal{T}}(k, 0, \mathbf{x}, \mathbf{x}') = \mathcal{T}_0 \delta(\mathbf{x} - \mathbf{x}'), \qquad \check{\mathcal{R}}(k, 0, \mathbf{x}, \mathbf{x}') = \mathcal{R}_0 \delta(\mathbf{x} - \mathbf{x}').$$

The symbol ○ stands for the Stratonovich stochastic integral in $z$, $B(z, \mathbf{x})$ is a real-valued Brownian field with covariance

$$\text{(3.5)} \qquad \mathbb{E}[B(z_1, \mathbf{x}_1)B(z_2, \mathbf{x}_2)] = \min\{z_1, z_2\} C_0(\mathbf{x}_1 - \mathbf{x}_2),$$

and we have used the notation

$$\text{(3.6)} \qquad C(z, \mathbf{x}) = \mathbb{E}[\nu(z' + z, \mathbf{x}' + \mathbf{x})\nu(z', \mathbf{x}')],$$

$$\text{(3.7)} \qquad C_0(\mathbf{x}) = \int_{-\infty}^{\infty} C(z, \mathbf{x}) \, dz.$$

The moments of the finite-dimensional distributions also converge, in the sense that

$$\text{(3.8)} \quad \begin{aligned} & \mathbb{E}\left[\prod_{j=1}^{q} p_{\text{tr}}^{\varepsilon}(s_j, \mathbf{x}_j)^{m_j} \prod_{j=1}^{\tilde{q}} p_{\text{ref}}^{\varepsilon}(\tilde{s}_j, \tilde{\mathbf{x}}_j)^{\tilde{m}_j}\right] \\ & \qquad \xrightarrow{\varepsilon \to 0} \mathbb{E}\left[\prod_{j=1}^{q} p_{\text{tr}}(s_j, \mathbf{x}_j)^{m_j} \prod_{j=1}^{\tilde{q}} p_{\text{ref}}(\tilde{s}_j, \tilde{\mathbf{x}}_j)^{\tilde{m}_j}\right], \end{aligned}$$

for any $q, \tilde{q} \in \mathbb{N}$, $(s_j)_{j=1,\ldots,q} \in \mathbb{R}^q$, $(\tilde{s}_j)_{j=1,\ldots,\tilde{q}} \in \mathbb{R}^{\tilde{q}}$, $(\mathbf{x}_j)_{j=1,\ldots,q} \in \mathbb{R}^{dq}$, $(\tilde{\mathbf{x}}_j)_{j=1,\ldots,\tilde{q}} \in \mathbb{R}^{d\tilde{q}}$, $(m_j)_{j=1,\ldots,q} \in \mathbb{N}^q$ and $(\tilde{m}_j)_{j=1,\ldots,\tilde{q}} \in \mathbb{N}^{\tilde{q}}$.

In [8] the existence and uniqueness have been established for the random process

$$V_k(z, \mathbf{x}) = \int \check{\mathcal{T}}(k, z, \mathbf{x}, \mathbf{x}') \phi(\mathbf{x}') \, d\mathbf{x}',$$

for a test function $\phi$ with unit $L^2(\mathbb{R}^d)$-norm, in the case $\mathcal{T}_0 = 1$. It is shown that the process $V_k(z, \mathbf{x})$ is a continuous Markov diffusion process on the unit ball of $L^2(\mathbb{R}^d, \mathbb{C})$. The moment equations moreover satisfy a closed system at each order [15]. The analysis can be readily extended to the pair $(p_{\text{tr}}, p_{\text{ref}})$ defined in terms of $(\check{\mathcal{T}}, \check{\mathcal{R}})$ and can be carried out jointly for all frequencies $k$ in an interval bounded away from 0 and infinity. We then get that the processes $p_{\text{tr}}$ and $p_{\text{ref}}$ have constant $L^2(\mathbb{R} \times \mathbb{R}^d)$-norms, so that the conservation of energy relation holds:

$$\text{(3.9)} \qquad \iint |p_{\text{tr}}(s, \mathbf{x})|^2 + |p_{\text{ref}}(s, \mathbf{x})|^2 \, ds \, d\mathbf{x} = \iint |b_{\text{inc}}(s, \mathbf{x})|^2 \, ds \, d\mathbf{x},$$

where we have also used the identity $\mathcal{R}_0^2 + \mathcal{T}_0^2 = 1$.

The main steps of the proof of Proposition 3.1 given in the next section are

(1) tightness and a priori estimates for $(p_{\text{tr}}^{\varepsilon}, p_{\text{ref}}^{\varepsilon})$,

(2) convergence of the finite-dimensional distributions of the process $(p_{\text{tr}}^{\varepsilon}, p_{\text{ref}}^{\varepsilon})$, using the convergence of specific moments of the reflection and transmission operators $\widehat{\mathcal{T}}^{\varepsilon}$ and $\widehat{\mathcal{R}}^{\varepsilon}$. Heuristically, the reasons for the convergence of $\widehat{\mathcal{T}}^{\varepsilon}$ and $\widehat{\mathcal{R}}^{\varepsilon}$ are that the terms with the rapid phases $\exp(\pm i 2kz/\varepsilon^4)$ in (2.20)–(2.21) vanish in the limit $\varepsilon \to 0$, and that the random "potential" in (2.22) can be replaced by a white noise. We then get formally the limit system (3.3)–(3.4). However, this holds true in a special weak sense only. Indeed, it is important to note that the reflection and transmission operators $\widehat{\mathcal{T}}^{\varepsilon}$ and $\widehat{\mathcal{R}}^{\varepsilon}$ themselves do not converge to $\widehat{\mathcal{T}}$ and $\widehat{\mathcal{R}}$ solution of (3.3)–(3.4), but only certain moments (expectations of products of components with distinct frequencies $k$), which are those needed to ensure the convergence of the transmitted and reflected fields. This approach is similar to the one used to prove the O'Doherty–Anstey theory in one-dimensional random media [7, 14] and to study the second-order statistics of the wave backscattered by a three-dimensional random medium [17, 18]. The limit moments are characterized by the system (A.15) of coupled equations, which gives the practical way to compute all the moments of the reflected and transmitted wave fields. Some particular applications will be given in Sections 4 and 5.

(3) use of Itô's lemma for Hilbert-space-valued processes [20], Theorem 2.4, in order to check that the specific moments of the reflection and transmission operators $\widehat{\mathcal{T}}$ and $\widehat{\mathcal{R}}$ given by (3.3)–(3.4) are solutions of the system (A.15) of coupled equations.

3.2. *Proof of Proposition 3.1.* This section is devoted to the proof of Proposition 3.1. We shall use a technique similar to the one presented in [14] in the case of randomly layered media.

*Step 1. A priori estimates.* From (2.12)–(2.13) we can check that, for any $k$, the integral

$$\int |\check{a}^{\varepsilon}(k,z,\mathbf{x})|^2 - |\check{b}^{\varepsilon}(k,z,\mathbf{x})|^2 \, d\mathbf{x}$$

is conserved in $z$. Applying this conservation relation at $z=0$ and $z=L$, and taking into account the boundary conditions (2.14)–(2.15), we obtain

$$\int |\check{a}^{\varepsilon}(k,L,\mathbf{x})|^2 \, d\mathbf{x} + (1-\mathcal{R}_0^2)\int |\check{b}^{\varepsilon}(k,0,\mathbf{x})|^2 \, d\mathbf{x} = \int |\check{b}_{\text{inc}}(k,\mathbf{x})|^2 \, d\mathbf{x}.$$

Using now (2.16)–(2.17) and the identity $\mathcal{R}_0^2 + \mathcal{T}_0^2 = 1$, we obtain

$$(3.10) \quad \int |\check{a}_1^{\varepsilon}(k,L,\mathbf{x})|^2 \, d\mathbf{x} + \int |\check{b}_0^{\varepsilon}(k,0,\mathbf{x})|^2 \, d\mathbf{x} = \int |\check{b}_{\text{inc}}(k,\mathbf{x})|^2 \, d\mathbf{x},$$

which expresses the fact that the power of the incoming wave is fully recovered by the transmitted and reflected waves. Integrating in $k$ and using

Parseval's equality establishes the total energy conservation relation

$$(3.11) \quad \iint |p_{\text{ref}}^\varepsilon(s,\mathbf{x})|^2 \, d\mathbf{x}\, ds + \iint |p_{\text{tr}}^\varepsilon(s,\mathbf{x})|^2 \, d\mathbf{x}\, ds = \iint |b_{\text{inc}}(s,\mathbf{x})|^2 \, d\mathbf{x}\, ds.$$

We first state a priori estimates for our quantities of interest.

LEMMA 3.1. *There exists $C > 0$ such that, uniformly in $\varepsilon$ and in $s_0, s_1$,*

$$(3.12) \quad \begin{aligned} \int |p_{\text{tr}}^\varepsilon(s_0,\mathbf{x})|^2 \, d\mathbf{x} &\leq C \quad \text{and} \\ \int |p_{\text{tr}}^\varepsilon(s_1,\mathbf{x}) - p_{\text{tr}}^\varepsilon(s_0,\mathbf{x})|^2 \, d\mathbf{x} &\leq C|s_1 - s_0|. \end{aligned}$$

*The same estimate holds true for $p_{\text{ref}}^\varepsilon$.*

PROOF. Using Sobolev's embedding $L^\infty(\mathbb{R}) \subset H^1(\mathbb{R})$, there exists a constant $C_{\text{sob}}$ such that, for any $\mathbf{x}$,

$$\sup_s |p_{\text{tr}}^\varepsilon(s,\mathbf{x})|^2 \leq C_{\text{sob}} \|p_{\text{tr}}^\varepsilon(\cdot,\mathbf{x})\|_{H^1(\mathbb{R},\mathbb{R})} = \frac{C_{\text{sob}}}{2\pi} \int (1+k^2)|\check{b}^\varepsilon(k,0,\mathbf{x})|^2 \, dk,$$

where we have also used Parseval's equality. Integrating in $\mathbf{x}$ and using the conservation equation (3.10) yields the first result of the lemma:

$$\sup_s \int |p_{\text{tr}}^\varepsilon(s,\mathbf{x})|^2 \, d\mathbf{x} \leq \int \sup_s |p_{\text{tr}}^\varepsilon(s,\mathbf{x})|^2 \, d\mathbf{x} \leq \frac{C_{\text{sob}}}{2\pi} \iint (1+k^2)|\check{b}_{\text{inc}}(k,\mathbf{x})|^2 \, d\mathbf{x}\, dk.$$

By the Cauchy–Schwarz inequality, we have

$$|p_{\text{tr}}^\varepsilon(s_1,\mathbf{x}) - p_{\text{tr}}^\varepsilon(s_0,\mathbf{x})|^2 = \left| \int_{s_0}^{s_1} \frac{\partial p_{\text{tr}}^\varepsilon}{\partial s}(s,\mathbf{x}) \, ds \right|^2 \leq \int_{s_0}^{s_1} ds \int_{s_0}^{s_1} \left| \frac{\partial p_{\text{tr}}^\varepsilon}{\partial s}(s,\mathbf{x}) \right|^2 ds$$

$$\leq |s_1 - s_0| \int \left| \frac{\partial p_{\text{tr}}^\varepsilon}{\partial s}(s,\mathbf{x}) \right|^2 ds.$$

The integral in $\mathbf{x}$ of the last term in the inequality can be bounded uniformly as above. The reflected field can be analyzed in the same way, which completes the proof. □

*Step 2.* The moments of the finite-dimensional distribution of $(p_{\text{tr}}^\varepsilon(s,\mathbf{x}), p_{\text{ref}}^\varepsilon(s,\mathbf{x}))$ converge to those of $(p_{\text{tr}}(s,\mathbf{x}), p_{\text{ref}}(s,\mathbf{x}))$. The general moment (3.8) of $p_{\text{tr}}^\varepsilon(s,\mathbf{x})$ can be expressed as the multiple integral

$$\mathbb{E}\left[ \prod_{j=1}^{q} p_{\text{tr}}^\varepsilon(s_j,\mathbf{x}_j)^{m_j} \prod_{j=1}^{\tilde{q}} p_{\text{ref}}^\varepsilon(\tilde{s}_j,\tilde{\mathbf{x}}_j)^{\tilde{m}_j} \right]$$

$$= \frac{1}{(2\pi)^{(N+M)(d+1)}}$$

$$\times \int \cdots \int \prod_{h=1}^{q} \prod_{j=1}^{m_h} d\boldsymbol{\kappa}'_{h,j} \, d\boldsymbol{\kappa}_{h,j} \, dk_{h,j}$$

$$\times \prod_{h=1}^{\tilde{q}} \prod_{j=1}^{\tilde{m}_h} d\tilde{\boldsymbol{\kappa}}'_{h,j} \, d\tilde{\boldsymbol{\kappa}}_{h,j} \, d\tilde{k}_{h,j}$$

$$\times \prod_{h,j} (\hat{b}_{\mathrm{inc}}(k_{h,j}, \boldsymbol{\kappa}'_{h,j}) e^{i(\boldsymbol{\kappa}_{h,j} \cdot \mathbf{x}_h - k_{h,j} s_h)})$$

$$\times \prod_{h,j} (\hat{b}_{\mathrm{inc}}(\tilde{k}_{h,j}, \tilde{\boldsymbol{\kappa}}'_{h,j}) e^{i(\tilde{\boldsymbol{\kappa}}_{h,j} \cdot \tilde{\mathbf{x}}_h - \tilde{k}_{h,j} \tilde{s}_h)})$$

$$\times \mathbb{E}\left[\prod_{h,j} \widehat{\boldsymbol{\mathcal{T}}}^{\varepsilon}(k_{h,j}, L, \boldsymbol{\kappa}_{h,j}, \boldsymbol{\kappa}'_{h,j}) \prod_{h,j} \widehat{\boldsymbol{\mathcal{R}}}^{\varepsilon}(\tilde{k}_{h,j}, L, \tilde{\boldsymbol{\kappa}}_{h,j}, \tilde{\boldsymbol{\kappa}}'_{h,j})\right],$$

for $N = \sum_{h=1}^{q} m_h$ and $M = \sum_{h=1}^{\tilde{q}} \tilde{m}_h$. Therefore, the convergence of the general moment of the transmitted and reflected wave fields in the white-noise limit will follow from the convergence of the following specific moments $\mathbb{E}[I^{\varepsilon}(L)]$ of the transmission and reflection operators, where

$$(3.13) \qquad I^{\varepsilon}(L) = \prod_{j=1}^{N} \widehat{\boldsymbol{\mathcal{T}}}^{\varepsilon}(k_j, L, \boldsymbol{\kappa}_j, \boldsymbol{\kappa}'_j) \prod_{j=1}^{M} \widehat{\boldsymbol{\mathcal{R}}}^{\varepsilon}(\tilde{k}_j, L, \tilde{\boldsymbol{\kappa}}_j, \tilde{\boldsymbol{\kappa}}'_j).$$

We call these moments "specific" because we restrict our attention to the case in which the frequencies $k_j$, $\tilde{k}_j$ are all distinct.

We use diffusion approximation theorems to obtain equations for the moments $I^{\varepsilon}$ in the limit $\varepsilon \to 0$. In the Appendix we show that

$$\lim_{\varepsilon \to 0} \mathbb{E}[I^{\varepsilon}(L)] = \mathbb{E}\left[\prod_{j=1}^{N} \widehat{\boldsymbol{\mathcal{T}}}(k_j, L, \boldsymbol{\kappa}_j, \boldsymbol{\kappa}'_j) \prod_{j=1}^{M} \widehat{\boldsymbol{\mathcal{R}}}(\tilde{k}_j, L, \tilde{\boldsymbol{\kappa}}_j, \tilde{\boldsymbol{\kappa}}'_j)\right],$$

when the right-hand side expectation is taken with respect to the following Itô–Schrödinger model for the transmission and reflection operators:

$$(3.14) \quad d\widehat{\boldsymbol{\mathcal{T}}}(k, z, \boldsymbol{\kappa}, \boldsymbol{\kappa}') = -\frac{i|\boldsymbol{\kappa}'|^2}{2k} \widehat{\boldsymbol{\mathcal{T}}}(k, z, \boldsymbol{\kappa}, \boldsymbol{\kappa}') \, dz - \frac{k^2 C_0(\mathbf{0})}{8} \widehat{\boldsymbol{\mathcal{T}}}(k, z, \boldsymbol{\kappa}, \boldsymbol{\kappa}') \, dz$$

$$+ \frac{ik}{2(2\pi)^d} \int \widehat{\boldsymbol{\mathcal{T}}}(k, z, \boldsymbol{\kappa}, \boldsymbol{\kappa}_1) \, d\hat{B}(z, \boldsymbol{\kappa}_1 - \boldsymbol{\kappa}') \, d\boldsymbol{\kappa}_1,$$

$$d\widehat{\boldsymbol{\mathcal{R}}}(k, z, \boldsymbol{\kappa}, \boldsymbol{\kappa}') = -\frac{i(|\boldsymbol{\kappa}|^2 + |\boldsymbol{\kappa}'|^2)}{2k} \widehat{\boldsymbol{\mathcal{R}}}(k, z, \boldsymbol{\kappa}, \boldsymbol{\kappa}') \, dz$$

$$- \frac{k^2 C_0(\mathbf{0})}{4} \widehat{\boldsymbol{\mathcal{R}}}(k, z, \boldsymbol{\kappa}, \boldsymbol{\kappa}') \, dz$$

$$
\begin{aligned}
(3.15) \qquad & -\frac{k^2}{4(2\pi)^d} \int \widehat{C}_0(\boldsymbol{\kappa}_1) \widehat{\mathcal{R}}(k,z,\boldsymbol{\kappa}-\boldsymbol{\kappa}_1, \boldsymbol{\kappa}'-\boldsymbol{\kappa}_1) \, d\boldsymbol{\kappa}_1 \, dz \\
& + \frac{ik}{2(2\pi)^d} \int \bigl( \widehat{\mathcal{R}}(k,z,\boldsymbol{\kappa},\boldsymbol{\kappa}_1) \, d\hat{B}(z,\boldsymbol{\kappa}_1-\boldsymbol{\kappa}') \\
& \qquad\qquad\qquad + \widehat{\mathcal{R}}(k,z,\boldsymbol{\kappa}_1,\boldsymbol{\kappa}') \, d\hat{B}(z,\boldsymbol{\kappa}-\boldsymbol{\kappa}_1) \bigr) \, d\boldsymbol{\kappa}_1,
\end{aligned}
$$

with the initial conditions $\widehat{\mathcal{T}}(k,0,\boldsymbol{\kappa},\boldsymbol{\kappa}') = \mathcal{T}_0 \delta(\boldsymbol{\kappa}-\boldsymbol{\kappa}')$ and $\widehat{\mathcal{R}}(k,0,\boldsymbol{\kappa},\boldsymbol{\kappa}') = \mathcal{R}_0 \delta(\boldsymbol{\kappa}-\boldsymbol{\kappa}')$. Here we have used the notations (3.6)–(3.7) and the Brownian field $\hat{B}$ has the following operator-valued spatial covariance:

$$(3.16) \quad \mathbb{E}[\hat{B}(z_1,\boldsymbol{\kappa}_1) \hat{B}(z_2,\boldsymbol{\kappa}_2)] = \min\{z_1,z_2\} (2\pi)^d \widehat{C}_0(\boldsymbol{\kappa}_1) \delta(\boldsymbol{\kappa}_1+\boldsymbol{\kappa}_2),$$

where

$$(3.17) \qquad \widehat{C}_0(\boldsymbol{\kappa}) = \int_{-\infty}^{\infty} \int_{\mathbb{R}^d} C(z,\mathbf{x}) e^{-i\boldsymbol{\kappa}\cdot\mathbf{x}} \, d\mathbf{x} \, dz.$$

The field $\hat{B}$ is the partial Fourier transform of the field $B$ defined in the statement of the proposition. Consider next the reflection operator in the original spatial variables:

$$(3.18) \quad \check{\mathcal{T}}(k,z,\mathbf{x},\mathbf{x}') = \frac{1}{(2\pi)^d} \iint e^{i(\boldsymbol{\kappa}\cdot\mathbf{x}-\boldsymbol{\kappa}'\cdot\mathbf{x}')} \widehat{\mathcal{T}}(k,z,\boldsymbol{\kappa},\boldsymbol{\kappa}') \, d\boldsymbol{\kappa} \, d\boldsymbol{\kappa}',$$

$$(3.19) \quad \check{\mathcal{R}}(k,z,\mathbf{x},\mathbf{x}') = \frac{1}{(2\pi)^d} \iint e^{i(\boldsymbol{\kappa}\cdot\mathbf{x}-\boldsymbol{\kappa}'\cdot\mathbf{x}')} \widehat{\mathcal{R}}(k,z,\boldsymbol{\kappa},\boldsymbol{\kappa}') \, d\boldsymbol{\kappa} \, d\boldsymbol{\kappa}'.$$

Then we find that this operator is weakly characterized by the following Itô–Schrödinger diffusion:

$$
\begin{aligned}
d\check{\mathcal{T}}(k,z,\mathbf{x},\mathbf{x}') ={}& \frac{i}{2k} \Delta_{\mathbf{x}'} \check{\mathcal{T}}(k,z,\mathbf{x},\mathbf{x}') \, dz - \frac{k^2 C_0(\mathbf{0})}{8} \check{\mathcal{T}}(k,z,\mathbf{x},\mathbf{x}') \, dz \\
& + \frac{ik}{2} \check{\mathcal{T}}(k,z,\mathbf{x},\mathbf{x}') \, dB(z,\mathbf{x}'), \\
d\check{\mathcal{R}}(k,z,\mathbf{x},\mathbf{x}') ={}& \frac{i}{2k} (\Delta_{\mathbf{x}} + \Delta_{\mathbf{x}'}) \check{\mathcal{R}}(k,z,\mathbf{x},\mathbf{x}') \, dz \\
& - \frac{k^2 (C_0(\mathbf{0}) + C_0(\mathbf{x}'-\mathbf{x}))}{4} \check{\mathcal{R}}(k,z,\mathbf{x},\mathbf{x}') \, dz \\
& + \frac{ik}{2} \check{\mathcal{R}}(k,z,\mathbf{x},\mathbf{x}') (dB(z,\mathbf{x}) + dB(z,\mathbf{x}')).
\end{aligned}
$$

In Stratonovich form this diffusion model becomes (3.3). This proves the last statement of the proposition (the convergence of the moments).

*Step 3. Convergence of* $(p_{\text{tr}}^\varepsilon, p_{\text{ref}}^\varepsilon)$ *to* $(p_{\text{tr}}, p_{\text{ref}})$ *in* $C^0(\mathbb{R}, L_w^2(\mathbb{R}^d, \mathbb{R}^2)) \cap L_w^2(\mathbb{R}, L_w^2(\mathbb{R}^d, \mathbb{R}^2))$. Here $L_w^2$ is the $L^2$ space equipped with the weak topology. Lemma 3.1 shows that the process $(p_{\text{tr}}^\varepsilon, p_{\text{ref}}^\varepsilon)$ is tight in $C^0(\mathbb{R}, L_w^2(\mathbb{R}^d, \mathbb{R}^2))$. Moreover, the first estimate in the lemma shows that, for any $L^2(\mathbb{R}^d, \mathbb{R})$-functions $\phi, \psi$, the random processes

$$X_\phi^\varepsilon(s) = \int p_{\text{tr}}^\varepsilon(s, \mathbf{x}) \phi(\mathbf{x}) \, d\mathbf{x},$$

$$Y_\psi^\varepsilon(s) = \int p_{\text{ref}}^\varepsilon(s, \mathbf{x}) \psi(\mathbf{x}) \, d\mathbf{x}$$

are uniformly bounded. Therefore, the finite-dimensional distributions are characterized by the moments of the form

$$\mathbb{E}\left[\prod_{j=1}^q X_{\phi_j}^\varepsilon(s_j)^{m_j} \prod_{j=1}^{\tilde{q}} Y_{\psi_j}^\varepsilon(\tilde{s}_j)^{\tilde{m}_j}\right],$$

where $q, \tilde{q} \in \mathbb{N}$, $m_j, \tilde{m}_j \in \mathbb{N}$, $s_j, \tilde{s}_j \in \mathbb{R}$, $\phi_j, \psi_j \in L^2(\mathbb{R}^d, \mathbb{R})$. These moments can be written as multiple integrals:

$$\mathbb{E}\left[\prod_{j=1}^q X_{\phi_j}^\varepsilon(s_j)^{m_j} \prod_{j=1}^{\tilde{q}} Y_{\psi_j}^\varepsilon(\tilde{s}_j)^{\tilde{m}_j}\right]$$

$$= \frac{1}{(2\pi)^{(N+M)(d+1)}}$$

$$\times \int \cdots \int \prod_{h=1}^q \prod_{j=1}^{m_h} d\boldsymbol{\kappa}_{h,j}' \, d\boldsymbol{\kappa}_{h,j} \, dk_{h,j} \times \prod_{h=1}^{\tilde{q}} \prod_{j=1}^{\tilde{m}_h} d\tilde{\boldsymbol{\kappa}}_{h,j}' \, d\tilde{\boldsymbol{\kappa}}_{h,j} \, d\tilde{k}_{h,j}$$

$$\times \prod_{h,j} (\hat{b}_{\text{inc}}(k_{h,j}, \boldsymbol{\kappa}_{h,j}') \overline{\hat{\phi}_h}(\boldsymbol{\kappa}_{h,j}) e^{-ik_{h,j} s_h})$$

$$\times \prod_{h,j} (\hat{b}_{\text{inc}}(\tilde{k}_{h,j}, \tilde{\boldsymbol{\kappa}}_{h,j}') \overline{\hat{\psi}_h}(\tilde{\boldsymbol{\kappa}}_{h,j}) e^{-i\tilde{k}_{h,j} \tilde{s}_h})$$

$$\times \mathbb{E}\left[\prod_{h,j} \widehat{\mathcal{T}}^\varepsilon(k_{h,j}, L, \boldsymbol{\kappa}_{h,j}, \boldsymbol{\kappa}_{h,j}') \prod_{h,j} \widehat{\mathcal{R}}^\varepsilon(\tilde{k}_{h,j}, L, \tilde{\boldsymbol{\kappa}}_{h,j}, \tilde{\boldsymbol{\kappa}}_{h,j}')\right],$$

for $N = \sum_{h=1}^q m_h$ and $M = \sum_{h=1}^{\tilde{q}} \tilde{m}_h$. Note that only specific moments of quantities of the form (3.13) appear (i.e., moments of products of the transmission and reflection operators at distinct $k$). The convergence of these specific moments therefore implies the convergence of the finite-dimensional distributions, hence the weak convergence in $C^0(\mathbb{R}, L_w^2(\mathbb{R}^d, \mathbb{R}^2))$. Furthermore, the estimate (3.11) shows that the processes are tight in $L_w^2(\mathbb{R}, L_w^2(\mathbb{R}^d, \mathbb{R}^2))$ (the unit ball is compact in the weak topology). This proves the weak convergence in $L_w^2(\mathbb{R}, L_w^2(\mathbb{R}^d, \mathbb{R}^2))$.

*Step 4. Convergence of $(p_{\mathrm{tr}}^\varepsilon, p_{\mathrm{ref}}^\varepsilon)$ to $(p_{\mathrm{tr}}, p_{\mathrm{ref}})$ in $C^0(\mathbb{R}, L^2(\mathbb{R}^d, \mathbb{R}^2)) \cap L^2(\mathbb{R}, L^2(\mathbb{R}^d, \mathbb{R}^2))$.* Here $L^2$ is the $L^2$ space equipped with the strong topology. Since the convergence has been proved in the weak topology, it is sufficient to show that the $L^2$-norm is preserved. On the one hand, the $L^2$-norms of the processes $(p_{\mathrm{tr}}^\varepsilon, p_{\mathrm{ref}}^\varepsilon)$ are deterministic, independent of $\varepsilon$, and given by (3.11). On the other hand, the limit process $(p_{\mathrm{tr}}, p_{\mathrm{ref}})$ has constant $L^2$-norm given by (3.9). This proves the first statement of the proposition and completes its proof.

**4. The Wigner distributions.** We first introduce the dimensionless autocorrelation function $\mathcal{C}$ of the fluctuations of the random medium

$$C(z, \mathbf{x}) = \sigma^2 \mathcal{C}\left(\frac{z}{l_z}, \frac{\mathbf{x}}{l_x}\right),$$

where $\sigma$ is the standard deviation of the fluctuations of the random medium and $l_z$ (resp. $l_x$) is the longitudinal (resp. transverse) correlation radius of the medium. With this representation we have

$$C_0(\mathbf{x}) = \sigma^2 l_z \mathcal{C}_0\left(\frac{\mathbf{x}}{l_x}\right), \qquad \widehat{C}_0(\mathbf{u}) = \sigma^2 l_z l_x^d \widehat{\mathcal{C}}_0(\mathbf{u} l_x).$$

We assume next that the power spectral density $\hat{\mathcal{C}}_0(\mathbf{u})$ decays fast enough so that $\int |\mathbf{u}|^2 \hat{\mathcal{C}}_0(\mathbf{u})\, d\mathbf{u}$ is finite. This means that the autocorrelation function $\mathcal{C}_0(\mathbf{x})$ is at least twice differentiable at $\mathbf{x} = \mathbf{0}$, which corresponds to a smooth random medium. For simplicity, we assume also that the random fluctuations are isotropic in the transverse directions, in the sense that the autocorrelation function $\mathcal{C}_0(\mathbf{x})$ depends only on $|\mathbf{x}|$.

4.1. *The Wigner distribution of the transmitted wave.* We now consider two frequencies $k_1$ and $k_2$ in a frequency band centered at $k$ and we define the two-frequency Wigner distribution of the transmission operator by

$$
\begin{aligned}
W^T_{k_1,k_2}&(z, \mathbf{x}, \mathbf{x}', \mathbf{q}, \mathbf{q}') \\
&= \iint e^{-i(\mathbf{q}\cdot\mathbf{y} + \mathbf{q}'\cdot\mathbf{y}')} \\
&\quad \times \mathbb{E}\Bigg[\check{\boldsymbol{\mathcal{T}}}\left(k_1, z, \frac{\sqrt{k}}{\sqrt{k_1}}\left(\mathbf{x} + \frac{\mathbf{y}}{2}\right), \frac{\sqrt{k}}{\sqrt{k_1}}\left(\mathbf{x}' + \frac{\mathbf{y}'}{2}\right)\right) \\
&\quad \times \overline{\boldsymbol{\mathcal{T}}}\left(k_2, z, \frac{\sqrt{k}}{\sqrt{k_2}}\left(\mathbf{x} - \frac{\mathbf{y}}{2}\right), \frac{\sqrt{k}}{\sqrt{k_2}}\left(\mathbf{x}' - \frac{\mathbf{y}'}{2}\right)\right)\Bigg] d\mathbf{y}\, d\mathbf{y}'.
\end{aligned}
\tag{4.1}
$$

Using the stochastic equation (3.3) and Itô's formula, we find that the Wigner distribution satisfies the closed system

$$\frac{\partial W^T_{k_1,k_2}}{\partial z} + \frac{\mathbf{q}'}{k} \cdot \nabla_{\mathbf{x}'} W^T_{k_1,k_2}$$

$$= -\frac{C_0(\mathbf{0})(k_1^2 + k_2^2)}{8} W_{k_1,k_2}^T$$
$$+ \frac{k_1 k_2}{4(2\pi)^d} \int \widehat{C}_0(\mathbf{u}) W_{k_1,k_2}^T\left(z, \mathbf{x}, \mathbf{x}', \mathbf{q}, \mathbf{q}' - \frac{1}{2}\left(\frac{\sqrt{k}}{\sqrt{k_1}} + \frac{\sqrt{k}}{\sqrt{k_2}}\right)\mathbf{u}\right)$$
$$\times e^{i\mathbf{u}\cdot\mathbf{x}'(\sqrt{k}/\sqrt{k_1} - \sqrt{k}/\sqrt{k_2})} d\mathbf{u},$$

starting from
$$W_{k_1,k_2}^T(z=0,\mathbf{x},\mathbf{x}',\mathbf{q},\mathbf{q}') = \mathcal{T}_0^2 (4\pi^2 k_1 k_2/k^2)^{d/2} \delta(\mathbf{x}-\mathbf{x}')\delta(\mathbf{q}+\mathbf{q}').$$

It is possible to solve this system and to find an integral representation of the two-frequency Wigner distribution by using the approach of [9]. However, we aim at focusing on spatial aspects in the next sections, and we shall simplify the algebra by assuming that the bandwidth $B$ of the incoming wave with carrier wavenumber $k_0$ is small. To describe this regime it is convenient to introduce

(4.2) $$\beta = \frac{\sigma^2 k_0^2 L l_z}{4}, \qquad \alpha = \frac{L}{k_0 l_x^2}, \qquad \alpha_0 = \frac{L}{k_0 r_0^2},$$

where $r_0$ is the initial beam width, $\beta$ describes the intensity of random scattering, while $\alpha$ and $\alpha_0$ represent the intensities of diffraction on respectively the scales of the medium variations and the input beam. Note that these parameters correspond to inverse Fresnel numbers (up to a factor $2\pi$) relative to respectively the lateral medium correlation scale and the aperture; below we shall consider explicitly the case with small medium Fresnel number. We assume that the bandwidth $B$ of the incoming wave is small in the sense that

(4.3) $$B \ll B_c, \qquad B_c := k_0 \min(1, \alpha^{-1}, \alpha_0^{-1}, \beta^{-1}).$$

In this regime we can approximate the two-frequency Wigner distribution by its behavior at the carrier frequency. That is, if $k_1, k_2$ lie in the spectrum of the incoming wave, the two-frequency Wigner distribution $W_{k_1,k_2}^T$ can be approximated by the simplified Wigner distribution $W^T$ that depends only on the carrier wavenumber $k_0$ and not on the lag $k_1 - k_2$. This Wigner distribution can be written in the form

$$W^T(z,\mathbf{x},\mathbf{x}',\mathbf{q},\mathbf{q}') = \mathcal{T}_0^2 (2\pi)^d \mathcal{W}^T\left(\frac{z}{L}, \frac{\mathbf{x}}{r_0}, \frac{\mathbf{x}'}{r_0}, \mathbf{q} l_x, \mathbf{q}' l_x\right),$$

where $\mathcal{W}^T$ satisfies

(4.4) $$\frac{\partial \mathcal{W}^T}{\partial \zeta} + \frac{\alpha l_x}{r_0} \mathbf{q}' \cdot \nabla_{\mathbf{x}'} \mathcal{W}^T = \frac{\beta}{(2\pi)^d} \int \hat{\mathcal{C}}_0(\mathbf{u})[\mathcal{W}^T(\mathbf{q}'-\mathbf{u}) - \mathcal{W}^T(\mathbf{q}')] d\mathbf{u},$$

starting from $\mathcal{W}^T(\zeta = 0, \mathbf{x}, \mathbf{x}', \mathbf{q}, \mathbf{q}') = \delta(\mathbf{x}-\mathbf{x}')\delta(\mathbf{q}+\mathbf{q}')$. We remark that these are in fact the classical equations of radiative transport for angularly resolved wave energy density [22]. In [22] the radiative transfer equations are written in the standard time-dependent form, while (4.4) is written in a time-harmonic form in a reference frame moving with the background velocity along the $z$-axis corresponding to the form set forth in [10].

In this context we have that $\beta \mathcal{C}_0(\mathbf{0})$ is the total scattering cross-section, $\beta \hat{\mathcal{C}}_0(\cdot)$ the differential scattering cross-section describing coupling of modes depending on their relative propagation directions and $(\alpha l_x/r_0)\mathbf{q}'$ a transport or velocity vector. It is clear from this equation that diffractive effects (characterized by the term $\mathbf{q}' \cdot \nabla_{\mathbf{x}'}$) are of order 1 if $\alpha l_x/r_0 \sim 1$ or equivalently $k_0 r_0 l_x \sim L$. In terms of Fresnel numbers, diffractive effects are of order 1 if

$$\alpha_e \equiv \frac{L}{k_0 a_e^2} \sim 1,$$

where we have defined the effective aperture $a_e$ by

$$a_e = \sqrt{l_x r_0}.$$

In Section 5 we shall see that a physically important and mathematically interesting regime corresponds to

$$\alpha_0 \ll \alpha_e \sim 1 \ll \alpha,$$

so that from the point of view of the "medium Fresnel" number we are in a Fraunhofer diffraction scaling, while from the point of view of the "source Fresnel" number we are in a Fresnel diffraction scaling and finally from the point of view of the effective Fresnel number, $1/(2\pi\alpha_e)$, we are in a general or scalar diffraction theory setup.

By taking a Fourier transform in $\mathbf{q}'$ and $\mathbf{x}'$, we obtain a transport equation that can be integrated and we find the following integral representation for $\mathcal{W}^T$:

$$\begin{aligned}(4.5)\quad &\mathcal{W}^T(\zeta, \mathbf{x}, \mathbf{x}', \mathbf{q}, \mathbf{q}') \\ &= \frac{\mathcal{T}_0^2}{(4\pi^2 \alpha)^d} \iint e^{-i(\mathbf{q}'+\mathbf{q})\cdot \boldsymbol{\eta}_1 - i(r_0(\mathbf{x}'-\mathbf{x})/(\alpha l_x)+\mathbf{q}\zeta)\cdot \boldsymbol{\eta}_2} \\ &\quad \times e^{\beta \int_0^\zeta \mathcal{C}_0(\boldsymbol{\eta}_1+\boldsymbol{\eta}_2\zeta')-\mathcal{C}_0(\mathbf{0})\, d\zeta'}\, d\boldsymbol{\eta}_1\, d\boldsymbol{\eta}_2.\end{aligned}$$

This expression will be used in Section 5.1 to compute and discuss the two-point statistics of the transmitted field.

4.2. *The Wigner distribution for the reflected wave.* We define the two-frequency Wigner distribution of the reflection operator by

$$W^R_{k_1,k_2}(z, \mathbf{x}, \mathbf{x}', \mathbf{q}, \mathbf{q}')$$

$$
\begin{aligned}
(4.6) \quad &= \iint e^{-i(\mathbf{q}\cdot\mathbf{y}+\mathbf{q}'\cdot\mathbf{y}')} \\
&\quad \times \mathbb{E}\bigg[\check{\mathcal{R}}\Big(k_1, z, \frac{\sqrt{k}}{\sqrt{k_1}}\Big(\mathbf{x}+\frac{\mathbf{y}}{2}\Big), \frac{\sqrt{k}}{\sqrt{k_1}}\Big(\mathbf{x}'+\frac{\mathbf{y}'}{2}\Big)\Big) \\
&\quad \times \overline{\check{\mathcal{R}}}\Big(k_2, z, \frac{\sqrt{k}}{\sqrt{k_2}}\Big(\mathbf{x}-\frac{\mathbf{y}}{2}\Big), \frac{\sqrt{k}}{\sqrt{k_2}}\Big(\mathbf{x}'-\frac{\mathbf{y}'}{2}\Big)\Big)\bigg]\,d\mathbf{y}\,d\mathbf{y}'.
\end{aligned}
$$

If the bandwidth of the incoming wave satisfies (4.3) and if $k_1, k_2$ lie in the spectrum of the wave, then we find by using (3.4) that the two-frequency Wigner distribution $W^R_{k_1,k_2}$ can be approximated by the simplified Wigner distribution $W^R$ that depends only on the carrier wavenumber $k_0$ and not on the lag $k_1 - k_2$. This Wigner distribution satisfies the closed system

$$
\frac{\partial W^R}{\partial z} + \frac{\mathbf{q}}{k_0}\cdot\nabla_\mathbf{x} W^R + \frac{\mathbf{q}'}{k_0}\cdot\nabla_{\mathbf{x}'} W^R
$$
$$
= \frac{k_0^2}{4(2\pi)^d}\int \widehat{C}_0(\mathbf{u})
$$
$$
\times \bigg[W^R(z,\mathbf{x},\mathbf{x}',\mathbf{q}-\mathbf{u},\mathbf{q}') + W^R(z,\mathbf{x},\mathbf{x}',\mathbf{q},\mathbf{q}'-\mathbf{u})
$$
$$
+ 2W^R\Big(z,\mathbf{x},\mathbf{x}',\mathbf{q}-\frac{1}{2}\mathbf{u},\mathbf{q}'-\frac{1}{2}\mathbf{u}\Big)\cos(\mathbf{u}\cdot(\mathbf{x}-\mathbf{x}'))
$$
$$
- 2W^R\Big(z,\mathbf{x},\mathbf{x}',\mathbf{q}-\frac{1}{2}\mathbf{u},\mathbf{q}'+\frac{1}{2}\mathbf{u}\Big)\cos(\mathbf{u}\cdot(\mathbf{x}-\mathbf{x}'))
$$
$$
- 2W^R(z,\mathbf{x},\mathbf{x}',\mathbf{q},\mathbf{q}')\bigg]d\mathbf{u},
$$

starting from $W^R(z=0,\mathbf{x},\mathbf{x}',\mathbf{q},\mathbf{q}') = \mathcal{R}_0^2(2\pi)^d\delta(\mathbf{x}-\mathbf{x}')\delta(\mathbf{q}+\mathbf{q}')$. We now cast the Wigner distribution in a suitable dimensionless form. We consider the following Fourier transform $V^R$ of the Wigner distribution $W^R$:

$$
W^R(z,\mathbf{x},\mathbf{x}',\mathbf{q},\mathbf{q}') = \frac{1}{(2\pi)^d}\int V^R\Big(z, \frac{\mathbf{q}+\mathbf{q}'}{2}, \mathbf{q}-\mathbf{q}', \mathbf{s}\Big)e^{i\mathbf{s}\cdot(\mathbf{x}'-\mathbf{x})}\,d\mathbf{s},
$$

which we introduce because the stationary maps that we will identify in Lemma 4.1, in the asymptotic regime $\alpha \to \infty$, have simple representations in this new frame. Note also that this ansatz incorporates the fact that $W^R$ does not depend on $\mathbf{x} + \mathbf{x}'$, only on $\mathbf{x} - \mathbf{x}'$, $\mathbf{q}$ and $\mathbf{q}'$, which follows from the stationarity of the random medium. The Fourier-transformed operator $V^R(z,\mathbf{q},\mathbf{r},\mathbf{s})$ has the form

$$
V^R(z,\mathbf{q},\mathbf{r},\mathbf{s}) = \mathcal{R}_0^2(\pi l_x)^d e^{iz\mathbf{r}\cdot\mathbf{s}/k_0}\mathcal{V}^R\Big(\frac{z}{L}, \mathbf{q}l_x, \mathbf{r}l_x, \mathbf{s}l_x\Big),
$$

where $\mathcal{V}^R$ is the solution of the dimensionless system

$$
\begin{aligned}
\frac{\partial \mathcal{V}^R}{\partial \zeta} &= \frac{\beta}{(2\pi)^d} \\
&\quad \times \int \hat{\mathcal{C}}_0(\mathbf{u}) \Big[ \mathcal{V}^R\Big(\zeta, \mathbf{q} - \frac{1}{2}\mathbf{u}, \mathbf{r} - \mathbf{u}, \mathbf{s}\Big) e^{-i\alpha \mathbf{s}\cdot\mathbf{u}\zeta} \\
&\qquad + \mathcal{V}^R\Big(\zeta, \mathbf{q} - \frac{1}{2}\mathbf{u}, \mathbf{r} + \mathbf{u}, \mathbf{s}\Big) e^{i\alpha \mathbf{s}\cdot\mathbf{u}\zeta} \\
&\qquad + \mathcal{V}^R\Big(\zeta, \mathbf{q} - \frac{1}{2}\mathbf{u}, \mathbf{r}, \mathbf{s} - \mathbf{u}\Big) e^{-i\alpha \mathbf{r}\cdot\mathbf{u}\zeta} \\
&\qquad + \mathcal{V}^R\Big(\zeta, \mathbf{q} - \frac{1}{2}\mathbf{u}, \mathbf{r}, \mathbf{s} + \mathbf{u}\Big) e^{i\alpha \mathbf{r}\cdot\mathbf{u}\zeta} - 2\mathcal{V}^R(\zeta, \mathbf{q}, \mathbf{r}, \mathbf{s}) \\
&\qquad - \mathcal{V}^R\Big(\zeta, \mathbf{q} - \frac{1}{2}\mathbf{u}, \mathbf{r} - \mathbf{u}, \mathbf{s} + \mathbf{u}\Big) e^{i\alpha[(\mathbf{r}-\mathbf{s})\cdot\mathbf{u} - |\mathbf{u}|^2]\zeta} \\
&\qquad - \mathcal{V}^R\Big(\zeta, \mathbf{q} - \frac{1}{2}\mathbf{u}, \mathbf{r} - \mathbf{u}, \mathbf{s} - \mathbf{u}\Big) e^{-i\alpha[(\mathbf{r}+\mathbf{s})\cdot\mathbf{u} + |\mathbf{u}|^2]\zeta} \Big] d\mathbf{u},
\end{aligned}
$$
(4.7)

starting from $\mathcal{V}^R(\zeta, \mathbf{q}, \mathbf{r}, \mathbf{s}) = \delta(\mathbf{q})$. The parameters $\alpha$ and $\beta$ are given by (4.2).

We want now to analyze the regime in which the transverse correlation length $l_x$ of the medium is smaller than the beam width $r_0$. More exactly, we assume from now on in this section that

(1) $r_0 \gg l_x$, which means that the transverse correlation length of the medium is small,

(2) $k_0 r_0 l_x \sim L$, which means that diffractive effects are of order 1.

These two conditions are equivalent to $\alpha_0 \ll \alpha_e \sim 1 \ll \alpha$. Note that in the previous section we established the fact that diffraction plays a role for the transmitted wave for a propagation distance $L$ of the order of $k_0 r_0 l_x$, which is smaller than the usual Rayleigh length $k_0 r_0^2$. This is a well-known result [19], and we shall deduce it for the reflected field in the analytic framework that we have set forth.

The rapid transverse variations regime is particularly interesting to study because $W^R$ has a multiscale behavior. In (4.7) this regime gives rise to rapid phases. The following proposition describes the asymptotic behavior of $\mathcal{V}^R$ as $\alpha \to \infty$. The presence of singular layers at $\mathbf{r} = \mathbf{0}$ and at $\mathbf{s} = \mathbf{0}$ requires particular attention and is responsible, for instance, for the enhanced backscattering phenomenon studied in Section 5.3. The situation with $\alpha$ large corresponds to a strong diffraction situation, at the scale of the lateral medium fluctuations. In general [part (1) in Lemma 4.1] the intensity of the reflection operator decays exponentially according to the parameter $\beta \mathcal{C}_0(0)$ corresponding to the total scattering cross section. This decay follows from a partial loss of coherence by random forward scattering. However, as articulated in parts (2) and (3) of the lemma below, the coupling of wave modes depends on the full medium autocorrelation function if we look at nearby specular reflection or small spatial offset frequencies. This coupling will be important when we analyze enhanced backscattering in Section 5.3.

LEMMA 4.1. *(1) For any* $\mathbf{r} \neq \mathbf{0}$, $\mathbf{s} \neq \mathbf{0}$:

$$\mathcal{V}^R(\zeta, \mathbf{q}, \mathbf{r}, \mathbf{s}) \stackrel{\alpha \to \infty}{\longrightarrow} \delta(\mathbf{q}) e^{-2\beta \mathcal{C}_0(\mathbf{0})\zeta}. \tag{4.8}$$

(2) *For any* $\mathbf{s} \neq \mathbf{0}$ *we have* $\mathcal{V}^R(\zeta, \mathbf{q}, \frac{\mathbf{r}}{\alpha}, \mathbf{s}) \stackrel{\alpha \to \infty}{\longrightarrow} \mathcal{V}^R_{\mathbf{r}}(\zeta, \mathbf{q})$ *where* $\mathcal{V}^R_{\mathbf{r}}(\zeta, \mathbf{q})$ *is solution of*

$$\frac{\partial \mathcal{V}^R_{\mathbf{r}}}{\partial \zeta} = \frac{2\beta}{(2\pi)^d} \int \hat{\mathcal{C}}_0(\mathbf{u}) \Big[ \mathcal{V}^R_{\mathbf{r}}\Big(\zeta, \mathbf{q} - \frac{1}{2}\mathbf{u}\Big) \cos(\mathbf{r} \cdot \mathbf{u}\zeta) - \mathcal{V}^R_{\mathbf{r}}(\zeta, \mathbf{q}) \Big] d\mathbf{u}, \tag{4.9}$$

*and is given explicitly by*

$$\mathcal{V}^R_{\mathbf{r}}(\zeta, \mathbf{q}) = \frac{1}{(2\pi)^d} \int e^{-i\mathbf{q}\cdot\mathbf{u}} e^{\beta \int_0^\zeta \mathcal{C}_0(\mathbf{u}/2+\mathbf{r}\zeta') + \mathcal{C}_0(\mathbf{u}/2-\mathbf{r}\zeta') - 2\mathcal{C}_0(\mathbf{0}) d\zeta'} d\mathbf{u}. \tag{4.10}$$

*Similarly, for any* $\mathbf{r} \neq \mathbf{0}$ *we have* $\mathcal{V}^R(\zeta, \mathbf{q}, \mathbf{r}, \frac{\mathbf{s}}{\alpha}) \stackrel{\alpha \to \infty}{\longrightarrow} \mathcal{V}^R_{\mathbf{s}}(\zeta, \mathbf{q})$.

(3) *For any* $\mathbf{r}$ *and* $\mathbf{s}$ *we have*

$$\mathcal{V}^R\Big(\zeta, \mathbf{q}, \frac{\mathbf{r}}{\alpha}, \frac{\mathbf{s}}{\alpha}\Big) \stackrel{\alpha \to \infty}{\longrightarrow} \mathcal{V}^R_{\mathbf{r}}(\zeta, \mathbf{q}) + \mathcal{V}^R_{\mathbf{s}}(\zeta, \mathbf{q}) - \delta(\mathbf{q}) e^{-2\beta\mathcal{C}_0(\mathbf{0})\zeta}. \tag{4.11}$$

PROOF. In case (1), the rapid phases cancel the contributions of all but the term $\mathcal{V}^R(\zeta, \mathbf{q}, \mathbf{r}, \mathbf{s})$ in (4.7), and we get

$$\frac{\partial \mathcal{V}^R}{\partial \zeta} = -2 \frac{\beta}{(2\pi)^d} \int \hat{\mathcal{C}}_0(\mathbf{u}) \mathcal{V}^R \, d\mathbf{u} = -2\beta \mathcal{C}_0(\mathbf{0}) \mathcal{V}^R,$$

which gives (4.8). In case (2), we obtain in the limit $\alpha \to \infty$ the simplified system

$$\frac{\partial \tilde{\mathcal{V}}^R_{\mathbf{r}}}{\partial \zeta} = \frac{\beta}{(2\pi)^d} \int \hat{\mathcal{C}}_0(\mathbf{u}) \Big[ \tilde{\mathcal{V}}^R_{\mathbf{r}}\Big(\zeta, \mathbf{q} - \frac{1}{2}\mathbf{u}, \mathbf{s} - \mathbf{u}\Big) e^{-i\mathbf{r}\cdot\mathbf{u}\zeta}$$
$$+ \tilde{\mathcal{V}}^R_{\mathbf{r}}\Big(\zeta, \mathbf{q} - \frac{1}{2}\mathbf{u}, \mathbf{s} + \mathbf{u}\Big) e^{i\mathbf{r}\cdot\mathbf{u}\zeta} - 2\tilde{\mathcal{V}}^R_{\mathbf{r}}(\zeta, \mathbf{q}, \mathbf{s}) \Big] d\mathbf{u}.$$

We then Fourier transform this equation in $\mathbf{q}$ and $\mathbf{s}$, and we obtain that the solution does not depend on $\mathbf{s}$, that it satisfies (4.9), and that it is given by (4.10).

In case (3), we obtain the simplified system for $\mathcal{V}^R_{\mathbf{r},\mathbf{s}}(\zeta,\mathbf{q}) = \lim_{\alpha\to\infty} \mathcal{V}^R(\zeta,\mathbf{q},\frac{\mathbf{r}}{\alpha},\frac{\mathbf{s}}{\alpha})$:

$$\frac{\partial \mathcal{V}^R_{\mathbf{r},\mathbf{s}}}{\partial \zeta} = \frac{2\beta}{(2\pi)^d} \int \hat{\mathcal{C}}_0(\mathbf{u}) \Big[ \mathcal{V}^R_{\mathbf{s}}\Big(\zeta, \mathbf{q}-\frac{1}{2}\mathbf{u}\Big) \cos(\mathbf{s}\cdot\mathbf{u}\zeta)$$
$$+ \mathcal{V}^R_{\mathbf{r}}\Big(\zeta, \mathbf{q}-\frac{1}{2}\mathbf{u}\Big) \cos(\mathbf{r}\cdot\mathbf{u}\zeta) - \mathcal{V}^R_{\mathbf{r},\mathbf{s}}(\zeta,\mathbf{q}) \Big] d\mathbf{u}.$$

Using (4.9) satisfied by $\mathcal{V}^R_{\mathbf{s}}$ and $\mathcal{V}^R_{\mathbf{r}}$, we get

$$\frac{\partial \mathcal{V}^R_{\mathbf{r},\mathbf{s}}}{\partial \zeta} = \frac{\partial \mathcal{V}^R_{\mathbf{r}}}{\partial \zeta} + \frac{\partial \mathcal{V}^R_{\mathbf{s}}}{\partial \zeta} + 2\beta\mathcal{C}_0(\mathbf{0})[\mathcal{V}^R_{\mathbf{r}} + \mathcal{V}^R_{\mathbf{s}} - \mathcal{V}^R_{\mathbf{r},\mathbf{s}}],$$

which yields (4.11). □

## 5. Two-point statistics of the transmitted and reflected fields.

5.1. *The transmitted field.* The results of Section 4.1 allow us to compute the two-point statistics of the transmitted field. We assume that

(a) the pulse has carrier frequency $k_0$ and it is narrowband in the sense that it satisfies (4.3),

(b) the input beam spatial profile is Gaussian with radius $r_0$,

$$(5.1) \qquad b_{\rm inc}(t,\mathbf{x}) = f_0(t) e^{-ik_0 t} \exp\Big(-\frac{|\mathbf{x}|^2}{r_0^2}\Big)$$

(suppressing the complex conjugate part here and below),

(c) the transverse correlation radius $l_x$ of the random fluctuations of the medium is much smaller than $r_0$ and $k_0 r_0 l_x \sim L$. As we have discussed after (4.4), this last condition ensures that diffractive effects are of order 1.

Under (a) and (b), we find that the autocorrelation function of the transmitted field defined by

$$A_{\rm tr}(s,t,\mathbf{x},\mathbf{y}) = \lim_{\varepsilon\to 0} \mathbb{E}\Big[ p^\varepsilon_{\rm tr}\Big(s+\frac{t}{2}, \mathbf{x}+\frac{\mathbf{y}}{2}\Big) p^\varepsilon_{\rm tr}\Big(s-\frac{t}{2}, \mathbf{x}-\frac{\mathbf{y}}{2}\Big) \Big]$$

has the form

$$A_{\rm tr}(s,t,\mathbf{x},\mathbf{y}) = \mathcal{T}_0^2 f_0\Big(s+\frac{t}{2}\Big) \overline{f_0}\Big(s-\frac{t}{2}\Big) e^{-ik_0 t} \Big(\frac{r_0^2}{8\pi}\Big)^{d/2}$$
$$(5.2) \qquad \times \int e^{-|\boldsymbol{\eta} L/k_0 + \mathbf{y}|^2/(2r_0^2) - r_0^2|\boldsymbol{\eta}|^2/8} e^{-i\boldsymbol{\eta}\cdot\mathbf{x}}$$
$$\times e^{k_0^2/4 \int_0^L \mathcal{C}_0(\boldsymbol{\eta} z/k_0 + \mathbf{y}) - \mathcal{C}_0(\mathbf{0})\, dz} \, d\boldsymbol{\eta}.$$

Under (a)–(c), we obtain

$$A_{\rm tr}(s,t,\mathbf{x},\mathbf{y}) = \mathcal{T}_0^2 f_0\Big(s+\frac{t}{2}\Big)\overline{f_0}\Big(s-\frac{t}{2}\Big) e^{-ik_0 t}\Big(\frac{r_0^2}{8\pi}\Big)^{d/2} \quad (5.3)$$
$$\times \int e^{-r_0^2|\boldsymbol{\eta}|^2/8 - |\mathbf{y}|^2/(2r_0^2)} e^{-i\boldsymbol{\eta}\cdot\mathbf{x}} e^{k_0^2/4 \int_0^L C_0(\boldsymbol{\eta} z/k_0+\mathbf{y})-C_0(\mathbf{0})\,dz}\, d\boldsymbol{\eta}.$$

If, moreover, random scattering is strong, in the sense that $\beta \gg 1$, or equivalently $k_0^2 C_0(\mathbf{0}) L \gg 1$, and if the autocorrelation function of the random fluctuations of the medium is twice differentiable at zero:

$$C_0(\mathbf{x}) = C_0(\mathbf{0}) - \frac{D}{2}|\mathbf{x}|^2 + o(|\mathbf{x}|^2),$$
(5.4)
$$D = -\frac{1}{d}\Delta C_0(\mathbf{0}) = -\frac{\sigma^2 l_z}{d l_x^2}\Delta\mathcal{C}_0(\mathbf{0}),$$

then we obtain that the autocorrelation function has the Gaussian shape

$$A_{\rm tr}(s,t,\mathbf{x},\mathbf{y}) = \mathcal{T}_0^2 f_0\Big(s+\frac{t}{2}\Big)\overline{f_0}\Big(s-\frac{t}{2}\Big)e^{-ik_0 t}\Big(\frac{r_0}{r_T(L)}\Big)^d \quad (5.5)$$
$$\times \exp\Big(-\frac{2|\mathbf{x}|^2}{r_T(L)^2} - \frac{|\mathbf{y}|^2}{2\rho_T(L)^2} + i\frac{\mathbf{x}\cdot\mathbf{y}}{\chi_T(L)^2}\Big).$$

The beam radius $r_T(L)$, the correlation radius $\rho_T(L)$ and the parameter $\chi_T(L)$ are characterized by

$$(5.6) \qquad r_T(L) = r_0\sqrt{1+\frac{DL^3}{3r_0^2}},$$

$$(5.7) \qquad \rho_T(L) = r_0 \frac{\sqrt{1+DL^3/(3r_0^2)}}{\sqrt{1+k_0^2 r_0^2 DL/4 + k_0^2 D^2 L^4/48}},$$

$$(5.8) \qquad \chi_T(L) = \frac{r_T(L)}{\sqrt{k_0 DL^2/2}},$$

where we have taken into account that $k_0 r_0^2 \gg L$ in the considered regime. Note in particular that the beam width increases at the anomalous rate $L^{3/2}$ (which was first obtained in the physical literature in [13] and confirmed mathematically in [11]). Furthermore, the lateral correlation radius decays to zero, which means that the beam becomes partially coherent. We finally remark that these results hold true in the case with a smooth random medium [with $C_0$ twice differentiable at $\mathbf{0}$ as in (5.4)]; the situation with a rough random medium will be addressed elsewhere.

### 5.2. The reflected field.
Under (a) and (b), the limit autocorrelation function of the reflected field defined by

$$A_{\text{ref}}(s,t,\mathbf{x},\mathbf{y}) = \lim_{\varepsilon \to 0} \mathbb{E}\left[p_{\text{ref}}^{\varepsilon}\left(s+\frac{t}{2},\mathbf{x}+\frac{\mathbf{y}}{2}\right)\overline{p_{\text{ref}}^{\varepsilon}\left(s-\frac{t}{2},\mathbf{x}-\frac{\mathbf{y}}{2}\right)}\right]$$

has the form

$$A_{\text{ref}}(s,t,\mathbf{x},\mathbf{y})$$
$$= \mathcal{R}_0^2 f_0\left(s+\frac{t}{2}\right)\overline{f_0}\left(s-\frac{t}{2}\right)e^{-ik_0 t}\left(\frac{r_0^2}{8\pi}\right)^{d/2}$$
(5.9)
$$\times \iiint d\boldsymbol{\eta}_1\,d\boldsymbol{\eta}_2\,d\boldsymbol{\eta}_3 e^{iL\boldsymbol{\eta}_3\cdot\boldsymbol{\eta}_2/k_0 - r_0^2(|\boldsymbol{\eta}_3|^2+|\boldsymbol{\eta}_2-2\boldsymbol{\eta}_1|^2)/8}$$
$$\times e^{-i\boldsymbol{\eta}_3\cdot\mathbf{x}+i\mathbf{y}\cdot(\boldsymbol{\eta}_1+\boldsymbol{\eta}_2/2)}\mathcal{V}^R(1,\boldsymbol{\eta}_1,\boldsymbol{\eta}_2,\boldsymbol{\eta}_3),$$

where the function $\mathcal{V}^R$ is the solution of the system (4.7). Under (a)–(c) the asymptotic behavior of the function $\mathcal{V}^R$ is determined by Lemma 4.1 and we find that the autocorrelation is given by

$$A_{\text{ref}}(s,t,\mathbf{x},\mathbf{y}) = \mathcal{R}_0^2 f_0\left(s+\frac{t}{2}\right)\overline{f_0}\left(s-\frac{t}{2}\right)e^{-ik_0 t}\left(\frac{r_0^2}{8\pi}\right)^{d/2}$$
(5.10)
$$\times \int e^{-r_0^2|\boldsymbol{\eta}|^2/8 - |\mathbf{y}|^2/(2r_0^2)}e^{-i\boldsymbol{\eta}\cdot\mathbf{x}}$$
$$\times e^{k_0^2/4 \int_0^{2L} C_0(\boldsymbol{\eta} z/k_0+\mathbf{y})-C_0(\mathbf{0})\,dz}\,d\boldsymbol{\eta}.$$

This is exactly the form (5.3) of the autocorrelation function of the transmitted wave, upon the substitution $2L$ for $L$. This shows that we would have obtained the same result if we had assumed that the backward propagation was independent of the forward propagation. The independent approach is valid in the regime in which the beam width is much larger than the transverse correlation radius of the fluctuations of the random medium, but it is not valid in other regimes, as can be seen by comparing the full expressions (5.2) and (5.9).

### 5.3. Enhanced backscattering.
The comparison of the autocorrelation function of the reflected wave and that of the transmitted wave for the propagation distance $2L$ shows that, in the regime $\alpha \gg 1$, there is no coherent effect building up between the forward and backward propagations. However, corrective terms show that there are some residual effects. In particular, we would like to show that the reflected intensity exhibits a singular picture in a very narrow cone, of angular width of order $\alpha^{-1}$, around the backscattered direction. This phenomenon, called enhanced backscattering or weak localization, is widely discussed in the physical literature [3, 27]. The physical observation is that, for an incoming quasi-monochromatic quasi-plane wave, the mean reflected power has a local maximum in the backscattered direction, which is twice as large as the mean reflected power in the other directions.

In this section, we assume that the incoming wave has the form

$$b_{\text{inc}}(t, \mathbf{x}) = f(t)e^{-ik_0 t}g_{\text{inc}}(\mathbf{x}),$$

that it is narrowband in the sense that it satisfies (4.3), and that it is nearly a plane wave, in the sense that $\hat{g}_{\text{inc}}(\boldsymbol{\kappa})$ is concentrated at some $\boldsymbol{\kappa}_{\text{inc}}$ (assumed to be different from the normal incident vector $\mathbf{0}$). By "concentrated" we mean that the angular width of the incoming beam is smaller than $\alpha^{-1}$. The reflected signal in the direction $\boldsymbol{\kappa}_0$ is

$$\check{p}^{\varepsilon}_{\text{ref}}(s, \boldsymbol{\kappa}_0) = \int p^{\varepsilon}_{\text{ref}}(s, \mathbf{x})e^{-i\boldsymbol{\kappa}_0 \cdot \mathbf{x}} d\mathbf{x}$$

$$= \frac{1}{2\pi} \int \widehat{\mathcal{R}^{\varepsilon}}(k, L, \boldsymbol{\kappa}_0, \boldsymbol{\kappa}')\hat{b}_{\text{inc}}(k, \boldsymbol{\kappa}')e^{-iks} dk.$$

The moment of the square modulus of $\check{p}^{\varepsilon}_{\text{ref}}(s, \boldsymbol{\kappa}_0)$ only involves specific moments of quantities of the form (3.13) (with distinct $k$). Therefore this moment converges to the one of the limit process $\check{p}_{\text{ref}}(s, \boldsymbol{\kappa}_0)$ defined as the Fourier transform in $\mathbf{x}$ of $p_{\text{ref}}(s, \mathbf{x})$ given by (3.2). This means that the mean reflected intensity in the direction $\boldsymbol{\kappa}_0$ converges to

$$\mathbb{E}[|\check{p}^{\varepsilon}_{\text{ref}}(s, \boldsymbol{\kappa}_0)|^2] \xrightarrow{\varepsilon \to 0} \mathcal{R}_0^2 |f(s)|^2 I^R(\boldsymbol{\kappa}_0),$$

$$I^R(\boldsymbol{\kappa}_0) = 2^{-d} l_x^d \int \mathcal{V}^R\left(1, \frac{\boldsymbol{\kappa}_0 - \boldsymbol{\kappa}_1'}{2} l_x, (\boldsymbol{\kappa}_0 + \boldsymbol{\kappa}_1')l_x, \mathbf{0}\right) |\hat{g}_{\text{inc}}(\boldsymbol{\kappa}_1')|^2 d\boldsymbol{\kappa}_1'.$$

Using the fact that $\hat{g}_{\text{inc}}(\boldsymbol{\kappa})$ is concentrated at $\boldsymbol{\kappa}_{\text{inc}}$, we get

$$(5.11) \qquad I^R(\boldsymbol{\kappa}_0) = P \mathcal{V}^R\left(1, \frac{\boldsymbol{\kappa}_0 - \boldsymbol{\kappa}_{\text{inc}}}{2} l_x, (\boldsymbol{\kappa}_0 + \boldsymbol{\kappa}_{\text{inc}})l_x, \mathbf{0}\right),$$

where $P = 2^{-d} l_x^d \int |\hat{g}_{\text{inc}}(\boldsymbol{\kappa}_1')|^2 d\boldsymbol{\kappa}_1'$. This formula gives the mean reflected intensity in the direction $\boldsymbol{\kappa}_0$ and is valid for arbitrary values of $\alpha$ and $\beta$. Let us consider the regime $\alpha \gg 1$. The mean reflected intensity far enough from the backscattered direction $-\boldsymbol{\kappa}_{\text{inc}}$ is of the form

$$I^R(\boldsymbol{\kappa}_0) = P\mathcal{V}_{\mathbf{0}}^R\left(1, \frac{\boldsymbol{\kappa}_0 - \boldsymbol{\kappa}_{\text{inc}}}{2} l_x\right) \qquad \text{for } |\boldsymbol{\kappa}_0 + \boldsymbol{\kappa}_{\text{inc}}|l_x \gg \alpha^{-1},$$

where we have used the second point of Lemma 4.1. In a narrow angular cone around the backscattered direction $-\boldsymbol{\kappa}_{\text{inc}}$, the reflected intensity is locally larger:

$$I^R(-\boldsymbol{\kappa}_{\text{inc}} + \alpha^{-1}\boldsymbol{\kappa}) = P[\mathcal{V}_{\mathbf{0}}^R(1, -\boldsymbol{\kappa}_{\text{inc}} l_x) + \mathcal{V}^R_{\boldsymbol{\kappa} l_x}(1, -\boldsymbol{\kappa}_{\text{inc}} l_x)],$$

where we have used the third point of Lemma 4.1. If we assume, additionally, that $\beta \gg 1$, then we have

$$I^R(\boldsymbol{\kappa}_0) = P(\pi \mathcal{D}\beta)^{-d/2} e^{-|\boldsymbol{\kappa}_0 - \boldsymbol{\kappa}_{\text{inc}}|^2 l_x^2/(4\mathcal{D}\beta)} \qquad \text{for } |\boldsymbol{\kappa}_0 + \boldsymbol{\kappa}_{\text{inc}}| l_x \gg \alpha^{-1},$$

where $\mathcal{D}$ is the dimensionless version of $D$ given by (5.4): $D = \sigma^2 l_z l_x^{-2} \mathcal{D}$. This formula gives the width of the diffusion cone around the specular direction $\boldsymbol{\kappa}_{\text{inc}}$:

$$(5.12) \qquad \Delta \boldsymbol{\kappa}_{\text{spec}} = \frac{2\sqrt{\mathcal{D}\beta}}{l_x} = \frac{\sqrt{\mathcal{D}}\sigma k_0 \sqrt{L l_z}}{l_x} = \sqrt{DL} k_0.$$

On the top of this broad cone, we have a narrow cone of relative maximum equal to 2 centered along the backscattered direction $-\boldsymbol{\kappa}_{\text{inc}}$:

$$I^R(-\boldsymbol{\kappa}_{\text{inc}} + \alpha^{-1}\boldsymbol{\kappa}) = P(\pi \mathcal{D}\beta)^{-d/2} e^{-|\boldsymbol{\kappa}_{\text{inc}}|^2 l_x^2/(\mathcal{D}\beta)} [1 + e^{-\mathcal{D}\beta|\boldsymbol{\kappa}|^2 l_x^2/3}].$$

This shows that the width of the enhanced backscattering cone is

$$(5.13) \qquad \Delta \boldsymbol{\kappa}_{\text{EBC}} = \frac{\sqrt{3}}{l_x \sqrt{\mathcal{D}\beta}\alpha} = \frac{2\sqrt{3} l_x}{\sqrt{\mathcal{D}}\sigma \sqrt{l_z L^3}} = \frac{2\sqrt{3}}{\sqrt{DL^3}}.$$

Note that the angular width $\Delta \theta_{\text{EBC}} = \Delta \boldsymbol{\kappa}_{\text{EBC}}/k_0$ of the cone is proportional to the wavelength, as predicted by physical arguments based on diagrammatic expansions [27].

## APPENDIX: DERIVATION OF THE LIMIT MOMENT EQUATIONS

The purpose of this appendix is to compute the limit of the expectation of (3.13) as $\varepsilon \to 0$, for distinct frequencies $k_j$, $\tilde{k}_j$. Using (2.20) and (2.21) we find

$$\frac{dI^\varepsilon}{dz}(z) = \sum_{j=1}^N \prod_{l=1, l\neq j}^N \widehat{\boldsymbol{\mathcal{T}}}^\varepsilon(k_l, z, \boldsymbol{\kappa}_l, \boldsymbol{\kappa}_l')$$

$$\times \prod_{l=1}^M \widehat{\boldsymbol{\mathcal{R}}}^\varepsilon(\tilde{k}_l, z, \tilde{\boldsymbol{\kappa}}_l, \tilde{\boldsymbol{\kappa}}_l')$$

$$\times \left\{ \int \widehat{\boldsymbol{\mathcal{T}}}^\varepsilon(k_j, z, \boldsymbol{\kappa}_j, \boldsymbol{\kappa}_a) \widehat{\boldsymbol{\mathcal{L}}}^\varepsilon(k_j, z, \boldsymbol{\kappa}_a, \boldsymbol{\kappa}_j') d\boldsymbol{\kappa}_a \right.$$

$$+ e^{2ik_j z/\varepsilon^4} \int \widehat{\boldsymbol{\mathcal{T}}}^\varepsilon(k_j, z, \boldsymbol{\kappa}_j, \boldsymbol{\kappa}_a) \widehat{\boldsymbol{\mathcal{L}}}^\varepsilon(k_j, z, \boldsymbol{\kappa}_a, \boldsymbol{\kappa}_b)$$

$$\left. \times \widehat{\boldsymbol{\mathcal{R}}}^\varepsilon(k_j, z, \boldsymbol{\kappa}_b, \boldsymbol{\kappa}_j') d\boldsymbol{\kappa}_a d\boldsymbol{\kappa}_b \right\}$$

$$+ \sum_{j=1}^M \prod_{l=1}^N \widehat{\boldsymbol{\mathcal{T}}}^\varepsilon(k_l, z, \boldsymbol{\kappa}_l, \boldsymbol{\kappa}_l')$$

(A.14)
$$\times \prod_{l=1\neq j}^{M} \widehat{\mathcal{R}}^{\varepsilon}(\tilde{k}_l, z, \tilde{\boldsymbol{\kappa}}_l, \tilde{\boldsymbol{\kappa}}'_l)$$
$$\times \left\{ e^{-2i\tilde{k}_j z/\varepsilon^4} \widehat{\mathcal{L}}^{\varepsilon}(\tilde{k}_j, z, \tilde{\boldsymbol{\kappa}}_j, \tilde{\boldsymbol{\kappa}}'_j) + e^{2i\tilde{k}_j z/\varepsilon^4} \right.$$
$$\times \iint \widehat{\mathcal{R}}^{\varepsilon}(\tilde{k}_j, z, \tilde{\boldsymbol{\kappa}}_j, \boldsymbol{\kappa}_a) \widehat{\mathcal{L}}^{\varepsilon}(\tilde{k}_j, z, \boldsymbol{\kappa}_a, \boldsymbol{\kappa}_b)$$
$$\times \widehat{\mathcal{R}}^{\varepsilon}(\tilde{k}_j, z, \boldsymbol{\kappa}_b, \tilde{\boldsymbol{\kappa}}'_j) \, d\boldsymbol{\kappa}_a \, d\boldsymbol{\kappa}_b$$
$$+ \int \widehat{\mathcal{L}}^{\varepsilon}(\tilde{k}_j, z, \tilde{\boldsymbol{\kappa}}_j, \boldsymbol{\kappa}_a) \widehat{\mathcal{R}}^{\varepsilon}(\tilde{k}_j, z, \boldsymbol{\kappa}_a, \tilde{\boldsymbol{\kappa}}'_j) \, d\boldsymbol{\kappa}_a$$
$$\left. + \int \widehat{\mathcal{R}}^{\varepsilon}(\tilde{k}_j, z, \tilde{\boldsymbol{\kappa}}_j, \boldsymbol{\kappa}_a) \widehat{\mathcal{L}}^{\varepsilon}(\tilde{k}_j, z, \boldsymbol{\kappa}_a, \tilde{\boldsymbol{\kappa}}'_j) \, d\boldsymbol{\kappa}_a \right\}.$$

We next apply the diffusion approximation to get limit equations for the moments; see [14] for background material on and related applications of the diffusion approximation. Observe that the random coefficients are rapidly fluctuating in view of (2.22). Those coefficients that are of order $\varepsilon^{-1}$ are centered and fluctuate on the scale $\varepsilon^2$; moreover they are assumed to be rapidly mixing, giving a white-noise scaling situation. Moreover, the rapid phase terms $\exp(\pm 2ikz/\varepsilon^4)$ lead to some cancellations between interacting terms. Here, the fact that the frequencies are distinct plays a key role. As a consequence, by applying diffusion approximation results, we obtain the equations for the moments $\mathbb{E}[I^{\varepsilon}]$ in the limit $\varepsilon \to 0$:

$$\bar{I}(z) = \lim_{\varepsilon \to 0} \mathbb{E}[I^{\varepsilon}(z)].$$

We obtain from (A.14) that $\bar{I}$ solves a system of integro-differential equations

$$\frac{d\bar{I}}{dz}(z) = -\frac{i}{2} \left( \sum_{j=1}^{N} \frac{|\boldsymbol{\kappa}'_j|^2}{k_j} + \sum_{j=1}^{M} \frac{|\tilde{\boldsymbol{\kappa}}_j|^2 + |\tilde{\boldsymbol{\kappa}}'_j|^2}{\tilde{k}_j} \right) \bar{I}(z)$$
$$- \frac{C_0(\mathbf{0})}{8} \left( \sum_{j=1}^{N} k_j^2 + 2 \sum_{j=1}^{M} \tilde{k}_j^2 \right) \bar{I}(z)$$
$$- \frac{1}{8(2\pi)^d} \int \widehat{C}_0(\boldsymbol{\kappa})$$
$$\times \left\{ \sum_{j=1}^{N} \sum_{l \neq j} k_j k_l \bar{I}(\boldsymbol{\kappa}'_j - \boldsymbol{\kappa}, \boldsymbol{\kappa}'_l + \boldsymbol{\kappa}) \right.$$

$$+ 2\sum_{j=1}^{M} \tilde{k}_j^2 \bar{I}(\widetilde{\boldsymbol{\kappa}}_j - \boldsymbol{\kappa}, \tilde{\boldsymbol{\kappa}}'_j - \boldsymbol{\kappa})$$

(A.15)
$$+ \sum_{j=1}^{M} \sum_{l \neq j} \tilde{k}_j \tilde{k}_l (\bar{I}(\tilde{\boldsymbol{\kappa}}_j - \boldsymbol{\kappa}, \tilde{\boldsymbol{\kappa}}'_l - \boldsymbol{\kappa})$$

$$+ \bar{I}(\tilde{\boldsymbol{\kappa}}_l - \boldsymbol{\kappa}, \tilde{\boldsymbol{\kappa}}'_j - \boldsymbol{\kappa})$$

$$+ \bar{I}(\tilde{\boldsymbol{\kappa}}_j - \boldsymbol{\kappa}, \tilde{\boldsymbol{\kappa}}_l + \boldsymbol{\kappa})$$

$$+ \bar{I}(\tilde{\boldsymbol{\kappa}}'_j - \boldsymbol{\kappa}, \tilde{\boldsymbol{\kappa}}'_l + \boldsymbol{\kappa}))$$

$$+ 2\sum_{j=1}^{N} \sum_{l=1}^{M} k_j \tilde{k}_l (\bar{I}(\boldsymbol{\kappa}'_j - \boldsymbol{\kappa}, \tilde{\boldsymbol{\kappa}}_l - \boldsymbol{\kappa})$$

$$+ \bar{I}(\boldsymbol{\kappa}'_j - \boldsymbol{\kappa}, \tilde{\boldsymbol{\kappa}}'_l + \boldsymbol{\kappa})) \bigg\} d\boldsymbol{\kappa},$$

where we only write the shifted arguments for $\bar{I}$. The initial conditions are $\bar{I}(k, \tilde{k}, \boldsymbol{\kappa}, \tilde{\boldsymbol{\kappa}}, \boldsymbol{\kappa}', \tilde{\boldsymbol{\kappa}}', z = 0) = \mathcal{T}_0^N \prod_{j=1}^{N} \delta(\boldsymbol{\kappa}_j - \boldsymbol{\kappa}'_j) \mathcal{R}_0^M \prod_{j=1}^{M} \delta(\tilde{\boldsymbol{\kappa}}_j - \tilde{\boldsymbol{\kappa}}'_j)$. Using in particular the relation

$$\mathbb{E}\left[ \iiint_0^{z_a} \int_0^{z_b} \lambda_a(s_a, \boldsymbol{\kappa}_a) \lambda_b(s_b, \boldsymbol{\kappa}_b) \, d\hat{B}(s_a, \boldsymbol{\kappa}_a) \, d\hat{B}(s_b, \boldsymbol{\kappa}_b) \, d\boldsymbol{\kappa}_a \, d\boldsymbol{\kappa}_b \right]$$

$$= \iint_0^{\min(z_a, z_b)} \mathbb{E}\left[ \lambda_a(s, \boldsymbol{\kappa}) \lambda_b(s, -\boldsymbol{\kappa}) \right] (2\pi)^d \widehat{C}_0(\boldsymbol{\kappa}) \, ds \, d\boldsymbol{\kappa},$$

we can then verify that

$$\bar{I}(z) = \mathbb{E}\left[ \prod_{j=1}^{N} \widehat{\mathcal{T}}(k_j, L, \boldsymbol{\kappa}_j, \boldsymbol{\kappa}'_j) \prod_{j=1}^{M} \widehat{\mathcal{R}}(\tilde{k}_j, L, \tilde{\boldsymbol{\kappa}}_j, \tilde{\boldsymbol{\kappa}}'_j) \right],$$

where the right-hand side expectation is taken with respect to the Itô–Schrödinger model for the transmission and reflection operators in (3.14)–(3.15).

LABORATOIRE DE PROBABILITÉS ET
  MODÈLES ALÉATOIRES
& LABORATOIRE JACQUES-LOUIS LIONS
UNIVERSITÉ PARIS VII
2 PLACE JUSSIEU
75251 PARIS CEDEX 05
FRANCE
E-MAIL: garnier@math.jussieu.fr

DEPARTMENT OF MATHEMATICS
UNIVERSITY OF CALIFORNIA
IRVINE, CALIFORNIA 92697
USA
E-MAIL: ksolna@math.uci.edu